# Linear Transceiver design for Downlink Multiuser MIMO Systems: Downlink-Interference Duality Approach

Tadilo Endeshaw Bogale, *Student Member, IEEE* and Luc Vandendorpe *Fellow, IEEE*

***Abstract*— This paper considers linear transceiver design for downlink multiuser multiple-input multiple-output (MIMO) systems. We examine different transceiver design problems. We focus on two groups of design problems. The first group is the weighted sum mean-square-error (WSMSE) (i.e., symbol-wise or user-wise WSMSE) minimization problems and the second group is the minimization of the maximum weighted mean-square-error (WMSE) (symbol-wise or user-wise WMSE) problems. The problems are examined for the practically relevant scenario where the power constraint is a combination of per base station (BS) antenna and per symbol (user), and the noise vector of each mobile station is a zero-mean circularly symmetric complex Gaussian random variable with arbitrary covariance matrix. For each of these problems, we propose a novel downlink-interference duality based iterative solution. Each of these problems is solved as follows. First, we establish a new mean-square-error (MSE) downlink-interference duality. Second, we formulate the power allocation part of the problem in the downlink channel as a Geometric Program (GP). Third, using the duality result and the solution of GP, we utilize alternating optimization technique to solve the original downlink problem. For the first group of problems, we have established symbol-wise and user-wise WSMSE downlink-interference duality. These duality are established by formulating the noise covariance matrices of the interference channels as fixed point functions. On the other hand, for the second group of problems, we have established symbol-wise and user-wise MSE downlink-interference duality. These duality are established by formulating the noise covariance matrices of the interference channels as marginally stable (convergent) discrete-time-switched systems. The proposed duality based iterative solutions can be extended straightforwardly to solve many other linear transceiver design problems. We also show that our MSE downlink-interference duality unify all existing MSE duality. In our simulation results, we have observed that the proposed duality based iterative algorithms utilize less total BS power than that of the existing algorithms.**

## I. Introduction

Multiple-input multiple-output (MIMO) is a promising technique to exploit the spectral efficiency of wireless channels. This spectral efficiency can be exploited by applying signal processing at the transmitter (precoder) and receiver (decoder). Signal processing is performed to meet a certain design criterion. It is well known that most practically relevant design problems such as weighted sum rate maximization, rate or signal-to-interference-plus-noise-ratio (SINR) balancing and rate or SINR constrained power minimization can be equivalently expressed as mean-square-error (MSE) based problems (see for example [1]). Because of this, the current paper examines MSE-based problems. In general, the uplink channel MSE-based problems are better understood than those of the downlink channel. Due to this fact, most literatures focus on solving the downlink MSE-based problems. The downlink MSE-based problems can be solved by direct approach as in [2], [3] or by uplink-downlink duality approach as in [4]–[6].

For a given downlink channel system model and its MSE-based problem, the idea behind uplink-downlink duality is first to create the virtual uplink channel by exchanging the roles of the transmitter and receiver, and then to enable the precoder/decoder transformation from uplink to downlink channel and vice versa by ensuring the same MSE in both channels. Once these two tasks are performed, the downlink MSE-based problems are examined as follows: When the global optimality of the dual uplink channel MSE-based problem is guaranteed, the duality approach simply transfers the optimal uplink channel precoder/decoder pairs from uplink to downlink channel (see for example the sum MSE minimization problem in [5]). When the global optimal solution of the dual uplink channel MSE-based problem can not ensured, the duality approach examines the downlink MSE-based problems by iteratively switching between the uplink and downlink channel problems (see for example the problems in [7]).

Several MSE-based problems have been examined by duality approach [4]–[8]. However, the duality of these papers are able to solve total BS power constrained MSE-based problems only. In a practical multi-antenna base station (BS) system, the maximum power of each BS antenna is limited [9]. In some scenario allocating different powers to different users (symbols) according to their priority or protection level has some interest. This motivates [10] to solve (robust) sum MSE-based problems with per antenna, user and symbol power constraints by duality approach. However, since the problems in [10] allocate the same MSE weight to all symbols (users), [10] ignores priority and fairness issues in terms of MSE.

In a multimedia communication, different types of information (for example, audio and video information) can be sent to a user (all users) simultaneously [11], [12]. In such a scenario, for successful transmission, more priority (power) could be given to symbols (users) corresponding to

The authors would like to thank BELSPO for the financial support of the IAP project BESTCOM in the framework of which this work has been achieved. Part of this work has been published in the ICASSP, Kyoto, Japan, Mar. 2012. Tadilo Endeshaw Bogale and Luc Vandendorpe are with the ICTEAM Institute, Université catholique de Louvain, Place du Levant 2, 1348 - Louvain La Neuve, Belgium. Email: {tadilo.bogale, luc.vandendorpe}@uclouvain.be, Phone: +3210478071, Fax: +3210472089.

the video information. Thus, for this scenario, the design criteria may incorporate fairness/priority and power constraints for each symbol (user). On the other hand, examining combined (per antenna and symbol (user)) power constraints may have practical interest (for example in network MIMO). For these reasons, the current paper generalizes the work of [10] by incorporating both symbol-wise and user-wise MSE fairness/priority, and combined (i.e., per antenna and symbol (user)) power constraints. We examine the following problems: Minimization of symbol-wise weighted sum mean-square-error (WSMSE) constrained with per BS antenna and symbol powers ($\mathcal{P}1$), minimization of user-wise WSMSE constrained with per BS antenna and user powers ($\mathcal{P}2$), minimization of the maximum symbol-wise weighted mean-square-error (WMSE) constrained with per BS antenna and symbol powers ($\mathcal{P}3$) and minimization of the maximum user-wise WMSE constrained with per BS antenna and user powers ($\mathcal{P}4$). Each of these problems is examined for the scenario where the noise vector of each mobile station (MS) is a zero-mean circularly symmetric complex Gaussian (ZMCSCG) random variable with arbitrary covariance matrix.

To the best of our knowledge, the problems $\mathcal{P}1$ - $\mathcal{P}4$ are non-convex. Furthermore, duality based solutions for these problems with our noise covariance matrix assumptions are not known. In the current paper, we propose duality based iterative solutions to solve the problems. Each of these problems is solved as follows. First, we establish a new MSE downlink-interference duality. Second, we formulate the power allocation part of the problem in the downlink channel as a Geometric Program (GP). Third, using the duality result and the solution of GP, we utilize alternating optimization technique to solve the original downlink problem. For the problems $\mathcal{P}1$ and $\mathcal{P}2$, the duality are established by formulating the noise covariance matrices of the interference channels as fixed point functions. For these two problems, the noise covariance matrices of the dual interference channels are computed by modifying the approach of [10] to $\mathcal{P}1$ and $\mathcal{P}2$ of the current paper. On the other hand, for the problems $\mathcal{P}3$ and $\mathcal{P}4$, the duality are established by formulating the noise covariance matrices of the interference channels as new marginally stable (convergent) discrete-time-switched systems. The proposed duality based iterative solutions can be extended straightforwardly to solve many other linear transceiver design problems. We also show that our MSE downlink-interference duality unify all existing MSE duality. In our simulation results, we have observed that the proposed duality based iterative algorithms utilize less total BS power than that of the existing algorithms. The main contributions of the current paper can thus be summarized as follows:

1) To solve the problems $\mathcal{P}1$ and $\mathcal{P}2$, we have established WSMSE downlink-interference duality by formulating the noise covariance matrices of the interference channels as fixed point functions. These noise covariance matrices are formulated by modifying the approach of [10] to $\mathcal{P}1$ and $\mathcal{P}2$ of the current paper. As will be clear later, for WSMSE-based problems with a total BS power constraint function, the proposed duality based algorithm requires less computation than that of the existing duality based algorithms.
2) To solve the problems $\mathcal{P}3$ and $\mathcal{P}4$, we have established novel MSE (symbol-wise and user-wise) downlink-interference duality by formulating the noise covariance matrices of the interference channels as marginally stable (convergent) discrete-time-switched systems. Furthermore, as will be shown later in Section IX, the proposed duality based iterative solutions can be extended straightforwardly to solve many other linear transceiver design problems. We also show that the MSE downlink-interference duality of the current paper is also applicable to solve total BS power based linear transceiver design problems. Thus, the current duality unify all existing MSE duality[1].
3) By employing the system model of [1] and [8], we formulate the power allocation parts of $\mathcal{P}1$ - $\mathcal{P}4$ as GPs. The GPs are formulated by applying the GP formulation approach of [1]. Consequently, we are able to solve our problems by alternating optimization technique [4], [7], [8], [10] (i.e., duality based iterative algorithm).
4) In our simulation results, we have observed that the proposed duality based iterative algorithms utilize less total BS power than that of the existing algorithms.

This paper is organized as follows. In Section II, multiuser MIMO downlink and virtual interference channel system models are presented. In Section III, we formulate our problems $\mathcal{P}1$ - $\mathcal{P}4$ and discuss the general framework of our duality based iterative solutions. Sections IV - VIII present the proposed duality based iterative solutions for solving these problems. The extension of our duality based iterative algorithms to other problems is discussed in Section IX. In Section X, computer simulations are used to compare the performance of the proposed duality algorithms with that of existing algorithms. Finally, conclusions are drawn in Section XI.

*Notations:* Upper/lower case boldface letters denote matrices/column vectors. The $\mathbf{X}_{(n,n)}$, $\mathbf{X}_{(n,:)}$, $\mathrm{tr}(\mathbf{X})$, $\mathbf{X}^T$, $\mathbf{X}^H$ and $\mathrm{E}(\mathbf{X})$ denote the $(n,n)$ element, $n$th row, trace, transpose, conjugate transpose and expected value of $\mathbf{X}$, respectively. $\mathbf{I}_n$ is an identity matrix of size $n \times n$ and $\mathbb{C}^{M\times M}$ ($\Re^{M\times M}$) represent spaces of $M \times M$ matrices with complex (real) entries. The diagonal and block-diagonal matrices are represented by $\mathrm{diag}(.)$ and $\mathrm{blkdiag}(.)$ respectively. Subject to is denoted by $\mathrm{s.t}$ and $(.)^\star$ denotes optimal solution. The superscripts $(.)^{DL}$ and $(.)^I$ denote downlink and interference, respectively.

## II. System Model

In this section, multiuser MIMO downlink and virtual interference channel system models are discussed which are shown in Fig. 1. In the downlink channel, the BS and $k$th MS are equipped with $N$ and $M_k$ antennas, respectively. The total number of MS antennas are thus $M = \sum_{k=1}^K M_k$. By denoting the symbol intended for the $k$th user as $\mathbf{d}_k \in \mathbb{C}^{S_k\times 1}$ and $S = \sum_{k=1}^K S_k$, the entire symbol can be written in a data vector $\mathbf{d} \in \mathbb{C}^{S\times 1}$ as $\mathbf{d} = [\mathbf{d}_1^T, \cdots, \mathbf{d}_K^T]^T$. The BS precodes $\mathbf{d}$ into an $N$ length vector by using its overall precoder matrix $\mathbf{B} = [\mathbf{b}_{11}, \cdots, \mathbf{b}_{KS_K}]$, where $\mathbf{b}_{ks} \in \mathbb{C}^{N\times 1}$ is the precoder vector of the BS for the $k$th MS $s$th symbol. The $k$th MS employs a receiver $\mathbf{w}_{ks}$ to estimate the symbol $d_{ks}$. We follow the same channel matrix notations as in [8]. The estimates of the $k$th MS $s$th symbol ($\hat{d}_{ks}$) and $k$th user ($\hat{\mathbf{d}}_k$) are given by

$$\hat{d}_{ks} = \mathbf{w}_{ks}^H(\mathbf{H}_k^H \sum_{i=1}^{K} \mathbf{B}_i \mathbf{d}_i + \mathbf{n}_k) = \mathbf{w}_{ks}^H(\mathbf{H}_k^H \mathbf{B}\mathbf{d} + \mathbf{n}_k) \quad (1)$$

$$\hat{\mathbf{d}}_k = \mathbf{W}_k^H(\mathbf{H}_k^H \mathbf{B}\mathbf{d} + \mathbf{n}_k) \quad (2)$$

where $\mathbf{H}_k^H \in \mathbb{C}^{M_k \times N}$ is the channel matrix between the BS and $k$th MS, $\mathbf{W}_k = [\mathbf{w}_{k1}, \cdots \mathbf{w}_{kS_k}]$, $\mathbf{B}_k = [\mathbf{b}_{k1}, \cdots \mathbf{b}_{kS_k}]$ and $\mathbf{n}_k$ is the $k$th MS additive noise. Without loss of generality, we can assume that the entries of $\mathbf{d}_k$ are independent and identically distributed (i.i.d) ZMCSCG random variables all with unit variance, i.e., $\mathrm{E}\{\mathbf{d}_k \mathbf{d}_k^H\} = \mathbf{I}_{S_k}$, $\mathrm{E}\{\mathbf{d}_k \mathbf{d}_i^H\} = \mathbf{0}$, $\forall i \neq k$, and $\mathrm{E}\{\mathbf{d}_k \mathbf{n}_i^H\} = \mathbf{0}$, $\forall i, k$. The $k$th MS noise vector is a ZMCSCG random variable with covariance matrix $\mathbf{R}_{nk} \in \mathbb{C}^{M_k \times M_k}$.

[1]Note that the existing MSE duality are established for a total BS power based linear transceiver design problems (see [1], [4], [5]).

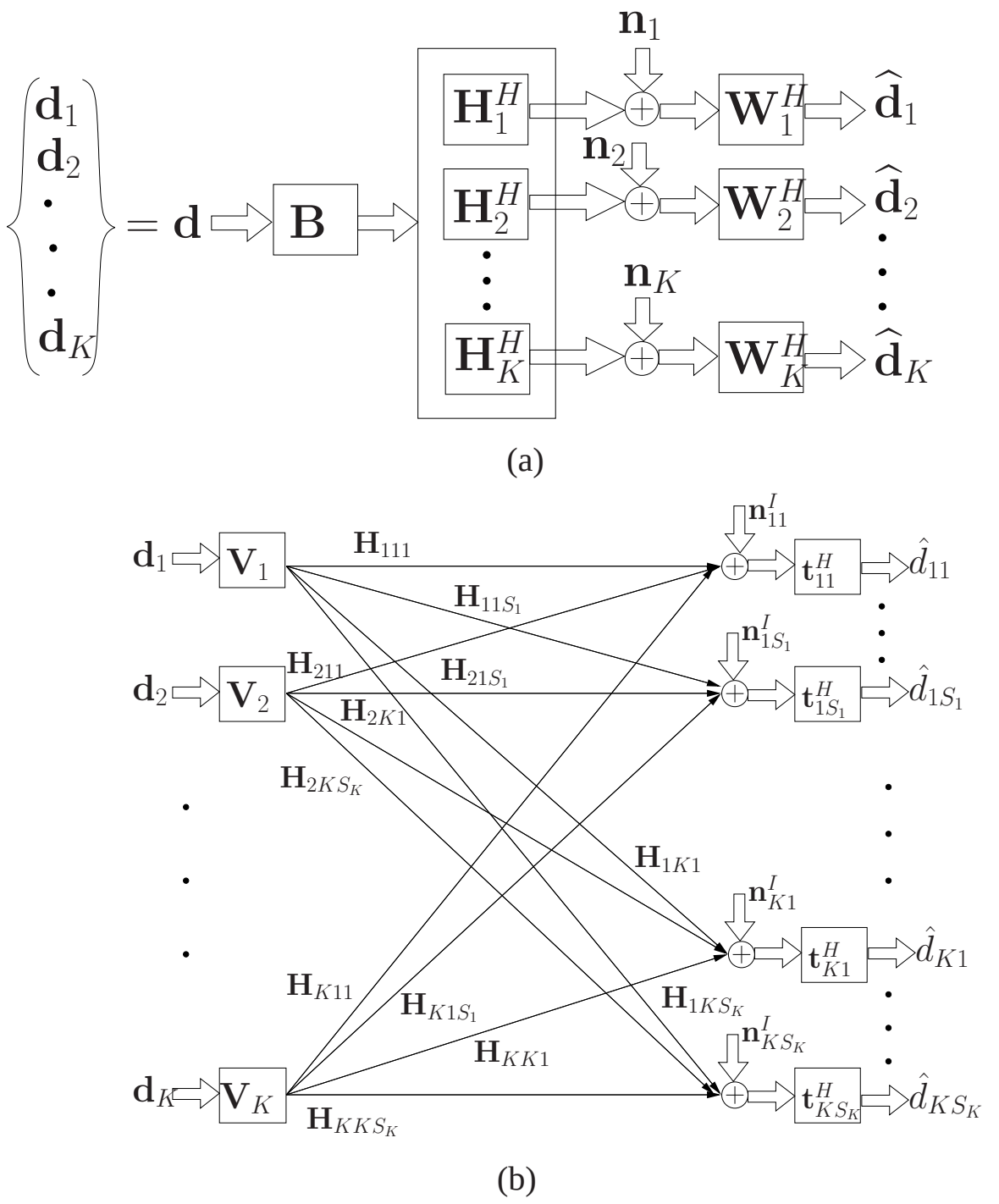

Fig. 1. Multiuser MIMO system model. (a) downlink channel. (b) virtual interference channel.

To establish our MSE downlink-interference duality, we model the virtual interference channel (Fig. 1.(b)) is modeled by introducing precoders $\{\mathbf{V}_k = [\mathbf{v}_{k1}, \cdots, \mathbf{v}_{kS_k}]\}_{k=1}^K$ and decoders $\{\mathbf{T}_k = [\mathbf{t}_{k1}, \cdots, \mathbf{t}_{kS_k}]\}_{k=1}^K$, where $\mathbf{v}_{ks} \in \mathbb{C}^{M_k \times 1}$ and $\mathbf{t}_{ks} \in \mathbb{C}^{N \times 1}, \forall k, s$. In this channel, it is assumed that the $k$th user's $s$th symbol ($d_{ks}$) is an i.i.d ZMCSCG random variable with variance $\zeta_{ks}$ and estimated independently by $\mathbf{t}_{ks} \in \mathcal{C}^{N\times 1}$, i.e., $\mathrm{E}\{d_{ks} d_{ks}^H\} = \zeta_{ks}$, $\mathrm{E}\{d_{ks} d_{ij}^H\} = 0$, $\forall (i,j) \neq (k,s)$, and $\mathrm{E}\{\mathbf{d}_k \mathbf{n}_i^H\} = \mathbf{0}$, $\forall i, k$. Moreover, $\{\mathbf{n}_{ks}^I, \forall s\}_{k=1}^K$ (Fig. 1.(b)) are also ZMCSCG random variables with covariance matrices $\{\mathbf{\Delta}_{ks} \in \Re^{N\times N} = \mathrm{diag}(\delta_{ks1}, \cdots, \delta_{ksN}), \forall s\}_{k=1}^K$ and the channels between the $k$th transmitter and all receivers are the same (i.e., $\{\mathbf{H}_{kjs} = \mathbf{H}_k, \forall j, s\}_{k=1}^K$). Note that from the system model aspect, the current paper and [10] share the same idea.

As can be seen from Fig. 1, the outputs of Fig. 1.(a) and Fig. 1.(b) are not the same. However, since Fig. 1.(b) is a "virtual" interference channel which is introduced just to solve the downlink MSE-based problems by duality approach, the output of Fig. 1.(b) is not required in practice. For this reason, the difference in the outputs of the downlink and interference channels of Fig. 1 will not affect the downlink MSE-based problem formulations and the duality based solutions.

For the downlink system model of Fig. 1.(a), the symbol-wise and user-wise MSEs can be expressed as

$$\begin{aligned}\xi_{ks}^{DL} =& \mathrm{E}_{\mathbf{d}}\{(\hat{d}_{ks} - d_{ks})(\hat{d}_{ks} - d_{ks})^H\} \\ =& \mathbf{w}_{ks}^H(\mathbf{H}_k^H \mathbf{B}\mathbf{B}^H \mathbf{H}_k + \mathbf{R}_{nk})\mathbf{w}_{ks} - \mathbf{w}_{ks}^H \mathbf{H}_k^H \mathbf{b}_{ks} - \\ & \mathbf{b}_{ks}^H \mathbf{H}_k \mathbf{w}_{ks} + 1 \end{aligned} \quad (3)$$

$$\begin{aligned}\xi_{k}^{DL} =& \mathrm{E}_{\mathbf{d}}\{(\hat{\mathbf{d}}_{k} - \mathbf{d}_{k})(\hat{\mathbf{d}}_{k} - \mathbf{d}_{k})^H\} \\ =& \mathrm{tr}\{\mathbf{I}_{S_k} + \mathbf{W}_k^H(\mathbf{H}_k^H \mathbf{B}\mathbf{B}^H \mathbf{H}_k + \mathbf{R}_{nk})\mathbf{W}_k - \\ & \mathbf{W}_k^H \mathbf{H}_k^H \mathbf{B}_k - \mathbf{B}_k^H \mathbf{H}_k \mathbf{W}_k\}. \end{aligned} \quad (4)$$

Using these two equations, the symbol-wise and user-wise WSMSEs can be expressed as

$$\begin{aligned}\xi_{ws}^{DL} =& \sum_{k=1}^{K}\sum_{s=1}^{S_k} \eta_{ks}\xi_{ks}^{DL} = \mathrm{tr}\{\boldsymbol{\eta} + \boldsymbol{\eta}\mathbf{W}^H\mathbf{H}^H\mathbf{B}\mathbf{B}^H\mathbf{H}\mathbf{W} + \\ & \boldsymbol{\eta}\mathbf{W}^H\mathbf{R}_n\mathbf{W} - \boldsymbol{\eta}\mathbf{W}^H\mathbf{H}^H\mathbf{B} - \boldsymbol{\eta}\mathbf{B}^H\mathbf{H}\mathbf{W}\} \end{aligned} \quad (5)$$

$$\begin{aligned}\xi_{wu}^{DL} =& \sum_{k=1}^{K} \tilde{\eta}_k \xi_k^{DL} = \mathrm{tr}\{\tilde{\boldsymbol{\eta}} + \tilde{\boldsymbol{\eta}}\mathbf{W}^H\mathbf{H}^H\mathbf{B}\mathbf{B}^H\mathbf{H}\mathbf{W} \\ & + \tilde{\boldsymbol{\eta}}\mathbf{W}^H\mathbf{R}_n\mathbf{W} - \tilde{\boldsymbol{\eta}}\mathbf{W}^H\mathbf{H}^H\mathbf{B} - \tilde{\boldsymbol{\eta}}\mathbf{B}^H\mathbf{H}\mathbf{W}\} \end{aligned} \quad (6)$$

where $\mathbf{R}_n = \mathrm{blkdiag}(\mathbf{R}_{n1}, \cdots, \mathbf{R}_{nK})$, $\boldsymbol{\eta} = \mathrm{diag}(\eta_{11}, \cdots, \eta_{1S_1}, \cdots, \eta_{K1}, \cdots, \eta_{KS_K})$ and $\tilde{\boldsymbol{\eta}} = \mathrm{blkdiag}(\tilde{\eta}_1 \mathbf{I}_{S_1}, \cdots, \tilde{\eta}_K \mathbf{I}_{S_K})$ with $\eta_{ks}$ and $\tilde{\eta}_k$ are the MSE weights of the $k$th user $s$th symbol and $k$th user, respectively. Like in the downlink channel, the interference channel symbol, user MSE and WSMSEs are expressed as

$$\begin{aligned}\xi_{ks}^{I} =& \mathbf{t}_{ks}^H \mathbf{\Gamma}_c \mathbf{t}_{ks} + \mathbf{t}_{ks}^H \mathbf{\Delta}_{ks}\mathbf{t}_{ks} - \mathbf{t}_{ks}^H \mathbf{H}_k \mathbf{v}_{ks}\zeta_{ks} - \\ & \zeta_{ks}\mathbf{v}_{ks}^H \mathbf{H}_k^H \mathbf{t}_{ks} + \zeta_{ks} \end{aligned} \quad (7)$$

$$\begin{aligned}\xi_{k}^{I} =& \mathrm{tr}\{\mathbf{T}_k^H \mathbf{\Gamma}_c \mathbf{T}_k - \mathbf{T}_k^H \mathbf{H}_k \mathbf{V}_k \boldsymbol{\zeta}_k - \boldsymbol{\zeta}_k \mathbf{V}_k^H \mathbf{H}_k^H \mathbf{T}_k + \boldsymbol{\zeta}_k\} + \\ & \sum_{s=1}^{S_k} \mathbf{t}_{ks}^H \mathbf{\Delta}_{ks} \mathbf{t}_{ks} \end{aligned} \quad (8)$$

$$\begin{aligned}\xi_{ws}^{I} =& \sum_{k=1}^{K}\sum_{s=1}^{S_k} \lambda_{ks}\xi_{ks}^{I} = \mathrm{tr}\{\boldsymbol{\lambda}\mathbf{T}^H\mathbf{\Gamma}_c\mathbf{T} - \boldsymbol{\lambda}\mathbf{T}^H\mathbf{H}\mathbf{V}\boldsymbol{\zeta} - \\ & \boldsymbol{\lambda}\boldsymbol{\zeta}\mathbf{V}^H\mathbf{H}^H\mathbf{T} + \boldsymbol{\lambda}\boldsymbol{\zeta}\} + \sum_{k=1}^{K}\sum_{s=1}^{S_k} \lambda_{ks}\mathbf{t}_{ks}^H\mathbf{\Delta}_{ks}\mathbf{t}_{ks} \end{aligned} \quad (9)$$

$$\begin{aligned}\xi_{wu}^{I} =& \sum_{k=1}^{K} \tilde{\lambda}_k \xi_k^{I} = \mathrm{tr}\{\tilde{\boldsymbol{\lambda}}\mathbf{T}^H\mathbf{\Gamma}_c\mathbf{T} - \tilde{\boldsymbol{\lambda}}\mathbf{T}^H\mathbf{H}\mathbf{V}\boldsymbol{\zeta} - \\ & \tilde{\boldsymbol{\lambda}}\boldsymbol{\zeta}\mathbf{V}^H\mathbf{H}^H\mathbf{T} + \tilde{\boldsymbol{\lambda}}\boldsymbol{\zeta}\} + \sum_{k=1}^{K}\sum_{s=1}^{S_k} \tilde{\lambda}_k\mathbf{t}_{ks}^H\mathbf{\Delta}_{ks}\mathbf{t}_{ks} \end{aligned} \quad (10)$$

where $\boldsymbol{\zeta}_k = \mathrm{diag}(\zeta_{k1}, \cdots, \zeta_{kS_k})$, $\boldsymbol{\zeta} = \mathrm{blkdiag}(\boldsymbol{\zeta}_1, \cdots, \boldsymbol{\zeta}_K)$, $\boldsymbol{\lambda} = \mathrm{diag}(\lambda_{11}, \cdots, \lambda_{1S_1}, \cdots, \lambda_{K1}, \cdots, \lambda_{KS_K})$, $\tilde{\boldsymbol{\lambda}} = \mathrm{blkdiag}(\tilde{\lambda}_1 \mathbf{I}_{S_1}, \cdots, \tilde{\lambda}_K \mathbf{I}_{S_K})$ and $\boldsymbol{\Gamma}_c = \sum_{i=1}^K \sum_{j=1}^{S_i} \zeta_{ij} \mathbf{H}_i \mathbf{v}_{ij} \mathbf{v}_{ij}^H \mathbf{H}_i^H$ with $\lambda_{ks}$ and $\tilde{\lambda}_k$ are the MSE weights of the $k$th user $s$th symbol and $k$th user, respectively.

## III. Problem Formulation

The aforementioned MSE-based optimization problems can be formulated as

$$\begin{aligned} \mathcal{P}1: & \min_{\{\mathbf{B}_k, \mathbf{W}_k\}_{k=1}^K} \sum_{k=1}^K \sum_{s=1}^{S_k} \eta_{ks} \xi_{ks}^{DL}, \\ & \text{s.t } [\mathbf{B}\mathbf{B}^H]_{(n,n)} \leq \check{p}_n, \ \mathbf{b}_{ks}^H \mathbf{b}_{ks} \leq \check{p}_{ks}, \ \forall n, k, s \end{aligned} \tag{11}$$

$$\begin{aligned} \mathcal{P}2: & \min_{\{\mathbf{B}_k, \mathbf{W}_k\}_{k=1}^K} \sum_{k=1}^K \tilde{\eta}_k \xi_k^{DL}, \\ & \text{s.t } [\mathbf{B}\mathbf{B}^H]_{(n,n)} \leq \check{p}_n, \ \mathrm{tr}\{\mathbf{B}_k^H \mathbf{B}_k\} \leq \hat{p}_k, \ \forall n, k \end{aligned} \tag{12}$$

$$\begin{aligned} \mathcal{P}3: & \min_{\{\mathbf{B}_k, \mathbf{W}_k\}_{k=1}^K} \max \ \rho_{ks} \xi_{ks}^{DL}, \\ & \text{s.t } [\mathbf{B}\mathbf{B}^H]_{(n,n)} \leq \check{p}_n, \ \mathbf{b}_{ks}^H \mathbf{b}_{ks} \leq \check{p}_{ks}, \ \forall n, k, s \end{aligned} \tag{13}$$

$$\begin{aligned} \mathcal{P}4: & \min_{\{\mathbf{B}_k, \mathbf{W}_k\}_{k=1}^K} \max \ \tilde{\rho}_k \xi_k^{DL}, \\ & \text{s.t } [\mathbf{B}\mathbf{B}^H]_{(n,n)} \leq \check{p}_n, \ \mathrm{tr}\{\mathbf{B}_k^H \mathbf{B}_k\} \leq \hat{p}_k, \ \forall n, k \end{aligned} \tag{14}$$

where $\tilde{\rho}_k(\hat{p}_k)$ and $\rho_{ks}(\check{p}_{ks})$ are the MSE balancing weights (maximum available power) of the $k$th user and $k$th user $s$th symbol, respectively, and $\check{p}_n$ denotes the maximum transmitted power by the $n$th antenna.

For both the WSMSE minimization and min max WMSE problems, different weights are given to different symbols (users). However, at optimality the solutions of these two problems are not necessarily the same. This is due to the fact that the aim of the WSMSE minimization problem is just to minimize the WSMSE of all symbols (users) (i.e., in such a problem the minimized WMSE of each symbol (user) depends on its corresponding channel gain), whereas the aim of min max WMSE problem is to minimize and balance the WMSE of each symbol (user) simultaneously (i.e., in such a problem all symbols (users) achieve the same minimized WMSE [13]). Moreover, as will be clear later, the solution approach of WSMSE minimization problem can not be extended straightforwardly to solve the min max WMSE problem. Due to these facts, we examine the WSMSE minimization and min max WMSE problems separately.

Since, the problems $\mathcal{P}1$ - $\mathcal{P}4$ are not convex, convex optimization framework can not be applied to solve them. To the best of our knowledge, duality based solutions for these problems are not known. In the following, we present an MSE downlink-interference duality based approach for solving each of these problems which is shown in **Algorithm I**[2].

**Algorithm I**

Initialization: For each problem, initialize $\{\mathbf{B}_k \neq \mathbf{0}\}_{k=1}^K$ such that the power constraint functions are satisfied[3]. Then, update $\{\mathbf{W}_k\}_{k=1}^K$ by using minimum mean-square-error (MMSE) receiver approach, i.e.,

$$\mathbf{W}_k = (\mathbf{H}_k^H \mathbf{B}\mathbf{B}^H \mathbf{H}_k + \mathbf{R}_{nk})^{-1} \mathbf{H}_k^H \mathbf{B}_k, \ \forall k. \tag{15}$$

**Repeat**: **Interference channel**

1) Transfer the symbol-wise (user-wise) WSMSE or WMSE from downlink to interference channel.
2) Update the receivers of the interference channel $\{\mathbf{t}_{ks}, \forall s\}_{k=1}^K$ using MMSE receiver technique.

**Downlink channel**

3) Transfer the symbol-wise (user-wise) WSMSE or WMSE from interference to downlink channel.
4) Update the receivers of the downlink channel $\{\mathbf{W}_k\}_{k=1}^K$ by MMSE receiver approach (15).

**Until** convergence.

The above iterative algorithm is already known in [5], [8] and [10]. However, the approaches of these papers can not ensure the power constraints of $\mathcal{P}1$ - $\mathcal{P}4$ at step 3 of **Algorithm I**. Hence, one can not apply the approaches of these papers to solve $\mathcal{P}1$ - $\mathcal{P}4$. In the following sections, we establish our MSE downlink-interference duality.

## IV. Symbol-wise WSMSE downlink-interference duality

This duality is established to solve symbol-wise WSMSE-based problems (for example $\mathcal{P}1$).

### A. Symbol-wise WSMSE transfer (From downlink to interference channel)

In order to use this WSMSE transfer for solving $\mathcal{P}1$, we set the interference channel precoder, decoder, noise covariance, input covariance and MSE weight matrices as

$$\mathbf{V} = \bar{\beta}\mathbf{W}, \ \mathbf{T} = \mathbf{B}/\bar{\beta}, \ \boldsymbol{\zeta} = \boldsymbol{\eta}, \ \boldsymbol{\lambda} = \mathbf{I}, \ \boldsymbol{\Delta}_{ks} = \boldsymbol{\Psi} + \mu_{ks}\mathbf{I} \tag{16}$$

where $\bar{\beta}$, $\{\psi_n\}_{n=1}^N$ and $\{\mu_{ks}, \forall s\}_{k=1}^K$ are positive real scalars that will be determined in the sequel and $\boldsymbol{\Psi} = \mathrm{diag}(\psi_1, \cdots, \psi_N)$. Substituting (16) into (9) and equating $\xi_{ws}^I = \xi_{ws}^{DL}$ yields

$$\begin{aligned} & \mathrm{tr}\{\mathbf{B}^H\mathbf{H}\mathbf{W}\boldsymbol{\eta}\mathbf{W}^H\mathbf{H}^H\mathbf{B} - \mathbf{B}^H\mathbf{H}\mathbf{W}\boldsymbol{\eta} - \boldsymbol{\eta}\mathbf{W}^H\mathbf{H}^H\mathbf{B} + \boldsymbol{\eta}\} + \\ & \frac{1}{\bar{\beta}^2} \sum_{k=1}^K \sum_{s=1}^{S_k} \mathbf{b}_{ks}^H (\boldsymbol{\Psi} + \mu_{ks}\mathbf{I}_N) \mathbf{b}_{ks} = \mathrm{tr}\{\boldsymbol{\eta}\mathbf{W}^H\mathbf{H}^H\mathbf{B}\mathbf{B}^H\mathbf{H}\mathbf{W} \\ & + \boldsymbol{\eta}\mathbf{W}^H\mathbf{R}_n\mathbf{W} - \boldsymbol{\eta}\mathbf{W}^H\mathbf{H}^H\mathbf{B} - \boldsymbol{\eta}\mathbf{B}^H\mathbf{H}\mathbf{W} + \boldsymbol{\eta}\}. \end{aligned}$$

It follows

$$\bar{\beta}^2 \tau = \sum_{n=1}^N \psi_n \check{p}_n + \sum_{k=1}^K \sum_{s=1}^{S_k} \mu_{ks} \bar{p}_{ks} = \check{\mathbf{p}}^T \boldsymbol{\psi} + \bar{\mathbf{p}}^T \boldsymbol{\mu} \tag{17}$$

where $\tau = \mathrm{tr}\{\boldsymbol{\eta}\mathbf{W}^H\mathbf{R}_n\mathbf{W}\}$, $\boldsymbol{\psi} = [\psi_1, \cdots, \psi_N]^T$, $\boldsymbol{\mu} = [\mu_{11}, \cdots, \mu_{1S_1}, \cdots, \mu_{K1}, \cdots, \mu_{KS_K}]^T$, $\check{\mathbf{p}} = [\check{p}_1, \cdots, \check{p}_N]^T$

[2] As will be clear later in Section VIII, to solve $\mathcal{P}3$ and $\mathcal{P}4$ (and more general MSE-based problems), an additional power allocation step is required. In **Algorithm I**, this step is omitted for clarity of presentation.

[3] For the simulation, we use $\{\mathbf{B}_k = [\mathbf{H}_k]_{(:,1:S_k)}\}_{k=1}^K$ followed by the appropriate normalization of $\{\mathbf{B}_k\}_{k=1}^K$ to ensure the power constraints.

and $\bar{\mathbf{p}} = [\bar{p}_{11}, \cdots, \bar{p}_{1S_1}, \cdots, \bar{p}_{K1}, \cdots, \bar{p}_{KS_K}]^T$ with $\bar{p}_{ks} = \mathbf{b}_{ks}^H \mathbf{b}_{ks}$, $\check{p}_n = \check{\mathbf{b}}_n^H \check{\mathbf{b}}_n$ and $\check{\mathbf{b}}_n^H$ is the $n$th row of $\mathbf{B}$.

The above equation shows that by choosing any $\{\psi_n\}_{n=1}^N$ and $\{\mu_{ks}, \forall s\}_{k=1}^K$ that satisfy (17), one can transfer the downlink channel precoder/decoder to the interference channel decoder/precoder ensuring $\xi_{ws}^{DL} = \xi_{ws}^{I_1}$, where $\xi_w^{I_1}$ is the interference WSMSE at step 1 of **Algorithm I**. However, here $\{\psi_n\}_{n=1}^N$ and $\{\mu_{ks}, \forall s\}_{k=1}^K$ should be selected in a way that $\mathcal{P}1$ can be solved by **Algorithm I**. To this end, we choose $\boldsymbol{\psi}$ and $\boldsymbol{\mu}$ as

$$\bar{\beta}^2 \tau \geq \check{\mathbf{p}}^T \boldsymbol{\psi} + \bar{\mathbf{p}}^T \boldsymbol{\mu}. \tag{18}$$

By doing so, the interference channel symbol-wise WSMSE is upper bounded by that of the downlink channel (i.e., $\xi_{ws}^{I_1} \leq \xi_{ws}^{DL}$). As will be clear later, to solve (11) with **Algorithm I**, $\bar{\beta}$, $\boldsymbol{\psi}$ and $\boldsymbol{\mu}$ should be selected as in (18). This shows that step 1 of **Algorithm I** can be carried out with (16). To perform step 2 of **Algorithm I**, we update $\mathbf{t}_{ks}$ of (16) by using the interference channel MMSE receiver approach which is expressed as

$$\begin{aligned} \mathbf{t}_{ks} =& (\boldsymbol{\Gamma}_c + \boldsymbol{\Delta}_{ks})^{-1} \mathbf{H}_k \mathbf{v}_{ks} \zeta_{ks} \\ =& \bar{\beta} (\mathbf{H} \mathbf{W} \boldsymbol{\eta} \mathbf{W}^H \mathbf{H}^H + \boldsymbol{\Psi} + \mu_{ks} \mathbf{I})^{-1} \mathbf{H}_k \mathbf{w}_{ks} \eta_{ks} \end{aligned} \tag{19}$$

where the second equality is obtained from (16). The above expression shows that by choosing $\{\mu_{ks} > 0, \forall s\}_{k=1}^K$, $\{\psi_n > 0\}_{n=1}^N$, we ensure $(\mathbf{H}\mathbf{W}\boldsymbol{\eta}\mathbf{W}^H\mathbf{H}^H + \boldsymbol{\Psi} + \mu_{ks}\mathbf{I})^{-1}$ exists. Next, we transfer the symbol-wise WSMSE from interference to downlink channel by ensuring the power constraint of $\mathcal{P}1$ (i.e., we perform step 3 of **Algorithm I**).

### *B. Symbol-wise WSMSE transfer (From interference to downlink channel)*

For a given symbol-wise WSMSE in the interference channel with $\boldsymbol{\zeta} = \boldsymbol{\eta}$ and $\boldsymbol{\lambda} = \mathbf{I}$, we can achieve the same WSMSE in the downlink channel (with the MSE weighting matrix $\boldsymbol{\eta}$) using a nonzero scaling factor ($\beta$) satisfying

$$\widetilde{\mathbf{B}} = \beta \mathbf{T}, \quad \widetilde{\mathbf{W}} = \mathbf{V}/\beta. \tag{20}$$

In this precoder/decoder transformation, we use the notations $\widetilde{\mathbf{B}}$ and $\widetilde{\mathbf{W}}$ to differentiate from the precoder and decoder matrices used in Section IV-A. By substituting (20) into $\xi_{ws}^{DL}$ (with $\widetilde{\mathbf{B}}$=$\mathbf{B}$, $\widetilde{\mathbf{W}}$=$\mathbf{W}$), equating the resulting symbol-wise WSMSE with that of the interference channel (9) and after some simple manipulations, we get

$$\begin{aligned} &\sum_{k=1}^K \sum_{s=1}^{S_k} \mathbf{t}_{ks}^H (\boldsymbol{\Psi} + \mu_{ks} \mathbf{I}_N) \mathbf{t}_{ks} = \frac{1}{\beta^2} \text{tr}\{\boldsymbol{\eta} \mathbf{V}^H \mathbf{R}_n \mathbf{V}\} \\ \Rightarrow \beta^2 =& \frac{\text{tr}\{\boldsymbol{\eta} \mathbf{V}^H \mathbf{R}_n \mathbf{V}\}}{\sum_{k=1}^K \sum_{s=1}^{S_k} \mathbf{t}_{ks}^H (\boldsymbol{\Psi} + \mu_{ks} \mathbf{I}_N) \mathbf{t}_{ks}} \\ =& \frac{\bar{\beta}^2 \tau}{\sum_{n=1}^N \psi_n \check{\mathbf{t}}_n^H \check{\mathbf{t}}_n + \sum_{k=1}^K \sum_{i=1}^{S_i} \mu_{ki} \mathbf{t}_{ki}^H \mathbf{t}_{ki}} \end{aligned} \tag{21}$$

where $\check{\mathbf{t}}_n^H$ is the $n$th row of the MMSE matrix $\mathbf{T}$ (19) and the third equality follows from (16). The power constraints of each BS antenna and symbol in the downlink channel are thus given by

$$\begin{aligned} \check{\widetilde{\mathbf{b}}}_n^H \check{\widetilde{\mathbf{b}}}_n =& \beta^2 \check{\mathbf{t}}_n^H \check{\mathbf{t}}_n \qquad (22) \\ =& \frac{\bar{\beta}^2 \tau \check{\mathbf{t}}_n^H \check{\mathbf{t}}_n}{\sum_{i=1}^N \psi_i \check{\mathbf{t}}_i^H \check{\mathbf{t}}_i + \sum_{i=1}^K \sum_{j=1}^{S_i} \mu_{ij} \mathbf{t}_{ij}^H \mathbf{t}_{ij}} \leq \check{p}_n, \ \forall n \\ \widetilde{\mathbf{b}}_{ks}^H \widetilde{\mathbf{b}}_{ks} =& \beta^2 \mathbf{t}_{ks}^H \mathbf{t}_{ks} \qquad (23) \\ =& \frac{\bar{\beta}^2 \tau \mathbf{t}_{ks}^H \mathbf{t}_{ks}}{\sum_{i=1}^N \psi_i \check{\mathbf{t}}_i^H \check{\mathbf{t}}_i + \sum_{i=1}^K \sum_{j=1}^{S_i} \mu_{ij} \mathbf{t}_{ij}^H \mathbf{t}_{ij}} \leq \check{p}_{ks}, \forall k, s \end{aligned}$$

where $\check{\widetilde{\mathbf{b}}}_n^H$ is the $n$th row of $\widetilde{\mathbf{B}}$. By multiplying both sides of (22) and (23) with $\psi_n, \forall n$ and $\mu_{ks}, \forall k, s$, we get

$$\psi_n \geq \check{f}_n \ \text{ and } \ \mu_{ks} \geq f_{ks}, \ \forall n, k, s \tag{24}$$

where $\check{f}_n = \frac{\bar{\beta}^2 \tau}{\check{p}_n} \frac{\psi_n \check{\mathbf{t}}_n^H \check{\mathbf{t}}_n}{\sum_{i=1}^N \psi_i \check{\mathbf{t}}_i^H \check{\mathbf{t}}_i + \sum_{i=1}^K \sum_{j=1}^{S_i} \mu_{ij} \mathbf{t}_{ij}^H \mathbf{t}_{ij}}$ and $f_{ks} = \frac{\bar{\beta}^2 \tau}{\check{p}_{ks}} \frac{\mu_{ks} \mathbf{t}_{ks}^H \mathbf{t}_{ks}}{\sum_{i=1}^N \psi_i \check{\mathbf{t}}_i^H \check{\mathbf{t}}_i + \sum_{i=1}^K \sum_{j=1}^{S_i} \mu_{ij} \mathbf{t}_{ij}^H \mathbf{t}_{ij}}$. Now, for any given $\bar{\beta}$, $\{\check{\mathbf{t}}_n^H \check{\mathbf{t}}_n\}_{n=1}^N$ and $\{\mathbf{t}_{ks}^H \mathbf{t}_{ks}, \forall s\}_{k=1}^K$, suppose that there exist $\{\psi_n > 0\}_{n=1}^N$ and $\{\mu_{ks} > 0, \forall s\}_{k=1}^K$ that satisfy

$$\psi_n = \check{f}_n \ \text{ and } \ \mu_{ks} = f_{ks}, \ \forall n, k, s. \tag{25}$$

From the above equation one can also achieve $\psi_n \check{p}_n = \check{f}_n \check{p}_n, \ \ \mu_{ks} \check{p}_{ks} = f_{ks} \check{p}_{ks} \ \forall n, k, s$. Summing up these expressions for all $n, k$ and $s$ results

$$\begin{aligned} \sum_{n=1}^N \psi_n \check{p}_n + \sum_{k=1}^K \sum_{s=1}^{S_k} \mu_{ks} \check{p}_{ks} =& \sum_{n=1}^N \check{f}_n \check{p}_n + \sum_{k=1}^K \sum_{s=1}^{S_k} f_{ks} \check{p}_{ks} \\ =& \bar{\beta}^2 \tau. \end{aligned} \tag{26}$$

This equation shows that the solution of (25) satisfies (26). Moreover, as $\{\check{p}_n \geq \bar{p}_n\}_{n=1}^N$ and $\{\check{p}_{ks} \geq \bar{p}_{ks}, \forall s\}_{k=1}^K$, the latter solution also ensures (18). Therefore, by choosing $\{\psi_n\}_{n=1}^N$ and $\{\mu_{ks}, \forall s\}_{k=1}^K$ such that (25) is satisfied, step 3 of **Algorithm I** can be performed. Furthermore, one can notice from (26) that $\bar{\beta}^2$ can be any positive value.

Next, we show that there exists at least a set of feasible $\{\psi_n > 0\}_{n=1}^N$ and $\{\mu_{ks} > 0, \forall s\}_{k=1}^K$ that satisfy (25). To this end, we consider the following Theorem [14].

*Theorem 1:* Let $(\mathbf{X}, \|.\|_2)$ be a complete metric space. We say that $\mathcal{F} : \mathbf{X} \rightarrow \mathbf{X}$ is an almost contraction, if there exist $\kappa(\tilde{\kappa}) \in [0, 1)$ and $\chi(\tilde{\chi}) \geq 0$ such that

$$\begin{aligned} &\|\mathcal{F}(\mathbf{x}) - \mathcal{F}(\mathbf{y})\|_2 \leq \kappa \|\mathbf{x} - \mathbf{y}\|_2 + \chi \|\mathbf{y} - \mathcal{F}(\mathbf{x})\|_2, \ \text{or} \qquad (27) \\ &\|\mathcal{F}(\mathbf{x}) - \mathcal{F}(\mathbf{y})\|_2 \leq \tilde{\kappa} \|\mathbf{x} - \mathbf{y}\|_2 + \tilde{\chi} \|\mathbf{x} - \mathcal{F}(\mathbf{y})\|_2, \ \forall \mathbf{x}, \mathbf{y} \in \mathbf{X}. \end{aligned}$$

If $\mathcal{F}$ satisfies (27), then the following holds true:

1) $\exists \mathbf{x} \in \mathbf{X} : \mathbf{x} = \mathcal{F}(\mathbf{x})$.
2) For any initial $\mathbf{x}_0 \in \mathbf{X}$, the iteration $\mathbf{x}_{n+1} = \mathcal{F}(\mathbf{x}_n)$ for n = 0, 1, 2, $\cdots$ converges to some $\mathbf{x}^\star \in \mathbf{X}$.
3) The solution $\mathbf{x}^\star$ is not necessarily unique.

*Proof:* See *Theorem 1.1* of [14]. Note that according to [15] (see (1.1) and (1.2) of [15]), the two inequalities of (27) are dual to each other.

Define $\mathbf{x}$ and $\mathcal{F}$ as $\mathbf{x} \triangleq [x_1, \cdots, x_{S+N}]^T = [\psi_1, \cdots, \psi_N, \mu_{11}, \cdots, \mu_{1S_1} \cdots, \mu_{K1}, \cdots, \mu_{KS_K}]^T$, $\mathcal{F}(\mathbf{x}) \triangleq [\check{f}_1, \cdots, \check{f}_N, f_{11}, \cdots, f_{1S_1}, \cdots, f_{K1}, \cdots, f_{KS_K}]$

with $\{x_n = \psi_n \in [\epsilon, (\bar{\beta}^2\tau - \epsilon\sum_{i=1,\ i\neq n}^{N} p_{im})/p_{nm}]\}_{n=1}^{N}$ and $\{x_r\}_{r=N+1}^{S+N} = \{\mu_{ks} =\in [\epsilon, (\bar{\beta}^2\tau - \epsilon\sum_{i=1}^{K}\sum_{j=1,(i,j)\neq(k,s)}^{S_i} p_{ijm})/p_{ksm}], \forall s\}_{k=1}^{K}$[4]. As we can see from (27), when $\|F(\mathbf{x}_1) - F(\mathbf{x}_2)\|_2 = 0$ with $\mathbf{x}_1 = \mathbf{x}_2$ or $\mathbf{x}_1 \neq \mathbf{x}_2$, one can set $\kappa(\tilde{\kappa}) = 0$ and $\chi(\tilde{\chi}) = 0$ to satisfy this inequality. And when $\|F(\mathbf{x}_1) - F(\mathbf{x}_2)\|_2 > 0$ (i.e., $\mathbf{x}_1 \neq \mathbf{x}_2$), one can select appropriate $\kappa(\tilde{\kappa}) \in [0, 1)$ and $\chi(\tilde{\chi}) \geq 0$ such that (27) is satisfied. This is due to the fact that in the latter case, $\|\mathbf{x}_2 - F(\mathbf{x}_1)\|_2 > 0$ and/or $\|\mathbf{x}_1 - F(\mathbf{x}_2)\|_2 > 0$ and $\|\mathbf{x}_1 - \mathbf{x}_2\|_2 > 0$ are positive and bounded. This explanation shows the existence of $\kappa(\tilde{\kappa}) \in [0, 1)$ and $\chi(\tilde{\chi}) \geq 0$ ensuring (27) for any $\|F(\mathbf{x}_1) - F(\mathbf{x}_2)\|_2$, $\mathbf{x}_1, \mathbf{x}_2 \in \mathbf{X}$ . Consequently, $F(\mathbf{x})$ is an almost contraction which implies

$$\mathbf{x}_{n+1} = F(\mathbf{x}_n),\ \mathbf{x}_0 = [x_{01}, x_{02}, \cdots, x_{0(S+N)}]^T \geq \epsilon\mathbf{1}_{N+S}, \quad \text{for } n = 0, 1, 2, \cdots \text{ converges} \tag{28}$$

where $\mathbf{1}_{N+S}$ is an $N + S$ length vector with each element equal to unity. Thus, there exist $\{\psi_n \geq \epsilon\}_{n=1}^{N}$ and $\{\mu_{ks} \geq \epsilon, \forall s\}_{k=1}^{K}$ that satisfy (25) and can be computed using (28). For numerical simulation we initialize $\mathbf{x}_0$ as $x_{01} = x_{02} = \cdots = x_{0(S+N)}$. However, finding the optimal initialization strategy is still an open research topic.

Once the appropriate $\{\psi_n\}_{n=1}^{N}$ and $\{\mu_{ks}\forall s\}_{k=1}^{K}$ are obtained, step 4 of **Algorithm I** is immediate and hence $\mathcal{P}1$ can be solved iteratively using this algorithm.

### C. Extension of the current duality for $\mathcal{P}1$ with a total BS power constraint

If the constraints of $\mathcal{P}1$ are modified to a total BS power, the power constraint at step 3 of **Algorithm I** can be ensured by applying the precoder/decoder transformation expression of [5]. The precoder/decoder transformation of [5] is performed by computing $S$ scaling factors. These scaling factors are obtained by solving $S$ systems of equations which require matrix inversion with complexity $O(S^3)$ (see (23) of [5]).

In the current paper, if the constraints of $\mathcal{P}1$ are modified to a total BS power, one can ensure the power constraint at step 3 of **Algorithm I** just by assigning $\Delta_{ks}$ of (16) as $\Delta_{ks} = \mathbf{I}$. By doing so, $\bar{\beta}^2$ of (17) and $\beta^2$ of (21) can be expressed as $\bar{\beta}^2 = \frac{\sum_{k=1}^{K}\sum_{s=1}^{S_k}\mathbf{b}_{ks}^H\mathbf{b}_{ks}}{\tau} = \frac{P_{max}}{\tau}$ and $\beta^2 = \frac{\text{tr}\{\boldsymbol{\eta}\mathbf{V}^H\mathbf{R}_n\mathbf{V}\}}{\sum_{k=1}^{K}\sum_{s=1}^{S_k}\mathbf{t}_{ks}^H\mathbf{t}_{ks}}$, where $P_{max}$ is the total BS power. Now by employing (20), the total BS power at step 3 of **Algorithm I** can thus be given as $\text{tr}\{\widetilde{\mathbf{B}}\widetilde{\mathbf{B}}^H\} = \beta^2\text{tr}\{\mathbf{T}\mathbf{T}^H\} = \bar{\beta}^2\tau = P_{max}$ (i.e., the total BS power constraint is satisfied). Thus, for $\mathcal{P}1$ (with a total BS power constraint), we do not need to use *Theorem I*. Moreover, our duality requires only one scaling factor to perform the precoder/decoder transformation (i.e., $\beta^2(\bar{\beta}^2)$). This shows that for this problem, the proposed duality based algorithm requires less computation compared to that of [5]. Note that the duality algorithm of [5] requires the same computation as that of [8] and less computation than that of [1] and [4]. Thus, it is sufficient to compare the current duality algorithm with the duality algorithm of [5].

[4]For our simulation, we use $\epsilon = \min(10^{-6}, \{\bar{\beta}\tau/p_{nm}\}_{n=1}^{N}, \{\bar{\beta}\tau/p_{ksm}, \forall s\}_{k=1}^{K})$.

For other WSMSE-based problems with a total BS power constraint function, the computational advantage of the current duality based algorithm over that of [5] can be analysed like in this subsection.

## V. User-wise WSMSE downlink-interference duality

This duality is established to solve user-wise WSMSE-based problems (for example $\mathcal{P}2$).

### A. User-wise WSMSE transfer (From downlink to interference channel)

To apply this WSMSE transfer for solving $\mathcal{P}2$, we set the precoder, decoder and noise covariance matrices as

$$\mathbf{V} = \tilde{\beta}\mathbf{W},\ \mathbf{T} = \mathbf{B}/\tilde{\beta},\ \boldsymbol{\zeta} = \tilde{\boldsymbol{\eta}},\ \tilde{\boldsymbol{\lambda}} = \mathbf{I}, \boldsymbol{\Delta}_{ks} = \boldsymbol{\Psi} + \mu_k\mathbf{I} \tag{29}$$

where $\tilde{\beta}$, $\{\psi_n\}_{n=1}^{N}$ and $\{\mu_k\}_{k=1}^{K}$ are real positive scalars. Substituting (29) into (10) and equating $\xi_{wu}^{I} = \xi_{wu}^{DL}$ yields

$$\tilde{\beta}^2\tilde{\tau} = \sum_{n=1}^{N}\psi_n\check{p}_n + \sum_{k=1}^{K}\mu_k p_k = \check{\mathbf{p}}^T\boldsymbol{\psi} + \tilde{\mathbf{p}}^T\tilde{\boldsymbol{\mu}} \tag{30}$$

where $\tilde{\tau} = \text{tr}\{\tilde{\boldsymbol{\eta}}\mathbf{W}^H\mathbf{R}_n\mathbf{W}\}$, $\tilde{\boldsymbol{\mu}} = [\mu_1, \cdots, \mu_K]^T$, $\tilde{\mathbf{p}} = [\tilde{p}_1, \cdots, \tilde{p}_K]^T$ with $\tilde{p}_k = \text{tr}\{\mathbf{B}_k\mathbf{B}_k^H\}$. Like in Section IV-A, we perform step 1 of **Algorithm I** by choosing $\tilde{\beta}^2$, $\boldsymbol{\psi}$ and $\tilde{\boldsymbol{\mu}}$ as

$$\tilde{\beta}\tilde{\tau} \geq \check{\mathbf{p}}^T\boldsymbol{\psi} + \tilde{\mathbf{p}}^T\tilde{\boldsymbol{\mu}}. \tag{31}$$

To perform step 2 of **Algorithm I**, we update $\mathbf{t}_{ks}$ of (29) using the interference channel MMSE receiver as

$$\mathbf{t}_{ks} = \tilde{\beta}(\mathbf{H}\mathbf{W}\tilde{\boldsymbol{\eta}}\mathbf{W}^H\mathbf{H}^H + \boldsymbol{\Psi} + \mu_k\mathbf{I})^{-1}\mathbf{H}_k\mathbf{w}_{ks}\tilde{\eta}_k. \tag{32}$$

This expression shows that by choosing $\{\mu_k > 0\}_{k=1}^{K}$, $\{\psi_n > 0\}_{n=1}^{N}$, we ensure that $(\mathbf{H}\mathbf{W}\boldsymbol{\eta}\mathbf{W}^H\mathbf{H}^H + \boldsymbol{\Psi} + \mu_k\mathbf{I})^{-1}$ exists.

### B. User-wise WSMSE transfer (From interference to downlink channel)

For a given user-wise WSMSE in the interference channel with $\boldsymbol{\zeta} = \tilde{\boldsymbol{\eta}}$ and $\tilde{\boldsymbol{\lambda}} = \mathbf{I}$, we can achieve the same WSMSE in the downlink channel (with the weighting matrix $\tilde{\boldsymbol{\eta}}$) by using a nonzero scaling factor ($\tilde{\tilde{\beta}}$) which satisfies

$$\widetilde{\mathbf{B}} = \tilde{\tilde{\beta}}\mathbf{T},\quad \widetilde{\mathbf{W}} = \mathbf{V}/\tilde{\tilde{\beta}}. \tag{33}$$

In this precoder/decoder transformation, we use the notations $\widetilde{\mathbf{B}}$ and $\widetilde{\mathbf{W}}$ to differentiate from the precoder and decoder matrices used in Section V-A. By substituting (33) into $\xi_{wu}^{DL}$ (with $\widetilde{\mathbf{B}}$=$\mathbf{B}$, $\widetilde{\mathbf{W}}$=$\mathbf{W}$), then equating the resulting user-wise WSMSE with that of the interference channel ($\xi_{wu}^{I}$) and after simple manipulations, we get

$$\tilde{\tilde{\beta}}^2 = \frac{\tilde{\beta}^2\tilde{\tau}}{\sum_{n=1}^{N}\psi_n\check{\mathbf{t}}_n^H\check{\mathbf{t}}_n + \sum_{k=1}^{K}\mu_k\text{tr}\{\mathbf{T}_k^H\mathbf{T}_k\}} \tag{34}$$

where $\check{\mathbf{t}}_n^H$ is the $n$th row of the MMSE matrix $\mathbf{T}$ (32). The power constraints of each BS antenna and user (i.e., step 3 of **Algorithm I**) in the downlink channel can be expressed as

$$\psi_n \geq \check{\tilde{f}}_n \quad \text{and} \quad \mu_k \geq \tilde{f}_k,\ \forall k \tag{35}$$

where

$$\check{\tilde{f}}_n = \frac{\tilde{\beta}^2\tilde{\tau}}{\breve{p}_n}\frac{\psi_n \check{\mathbf{t}}_n^H\check{\mathbf{t}}_n}{\sum_{i=1}^N \psi_i \check{\mathbf{t}}_i^H\check{\mathbf{t}}_i + \sum_{i=1}^K \mu_i \mathrm{tr}\{\mathbf{T}_i^H\mathbf{T}_i\}} \tag{36}$$

$$\tilde{f}_k = \frac{\tilde{\beta}^2\tilde{\tau}}{\hat{p}_k}\frac{\mu_k \mathrm{tr}\{\mathbf{T}_k^H\mathbf{T}_k\}}{\sum_{i=1}^N \psi_i \check{\mathbf{t}}_i^H\check{\mathbf{t}}_i + \sum_{i=1}^K \mu_i \mathrm{tr}\{\mathbf{T}_i^H\mathbf{T}_i\}}. \tag{37}$$

For given $\tilde{\beta}$, $\{\check{\mathbf{t}}_n^H\check{\mathbf{t}}_n\}_{n=1}^N$ and $\{\mathrm{tr}\{\mathbf{T}_k^H\mathbf{T}_k\}\}_{k=1}^K$, one can show that there exist $\{\psi_n\}_{n=1}^N$ and $\{\mu_k\}_{k=1}^K$ which satisfy

$$\psi_n = \check{\tilde{f}}_n \ \text{ and } \ \mu_k = \tilde{f}_k, \ \forall n,k. \tag{38}$$

The solution of (38) can be obtained exactly like that of (25). As $\{\breve{p}_n \geq \check{p}_n\}_{n=1}^N$ and $\{\hat{p}_k \geq \tilde{p}_k\}_{k=1}^K$, the latter solution also satisfies (31). Thus, $\mathcal{P}2$ can be solved using **Algorithm I**.

## VI. Symbol-wise MSE downlink-interference duality

In this section, we establish the symbol-wise MSE duality between downlink and interference channels. If all symbols are active, this duality can be applied to solve MSE based problems. However, as will be clear later, this duality requires more computation compared to the duality of Sections IV and V. Thus, we propose this duality to be employed for problems like in $\mathcal{P}3$ since this problem maintains all symbols active and can not be solved by the duality in Sections IV and V.

### *A. Symbol-wise MSE transfer (From downlink to interference channel)*

To apply this duality for $\mathcal{P}3$, we set the interference channel precoder, decoder, noise covariance, input covariance and MSE weight matrices as

$$\mathbf{v}_{ks} = \bar{\beta}_{ks}\mathbf{w}_{ks}, \ \mathbf{t}_{ks} = \mathbf{b}_{ks}/\bar{\beta}_{ks}, \boldsymbol{\zeta} = \mathbf{I}, \\ \boldsymbol{\Delta}_{ks} = \boldsymbol{\Psi} + \mu_{ks}\mathbf{I}_N, \forall k,s. \tag{39}$$

Substituting (39) into (7) and $\{\xi_{ks}^{DL} = \xi_{ks}^{I}, \forall s\}_{k=1}^K$ yields

$$\mathbf{w}_{ks}^H(\mathbf{H}_k^H\sum_{i=1}^K\sum_{j=1}^{S_i}\mathbf{b}_{i,j}\mathbf{b}_{i,j}^H\mathbf{H}_k + \mathbf{R}_{nk})\mathbf{w}_{ks} - \mathbf{w}_{ks}^H\mathbf{H}_k^H\mathbf{b}_{ks} \\ - \mathbf{b}_{ks}^H\mathbf{H}_k\mathbf{w}_{ks} + 1 = \frac{1}{\bar{\beta}_{ks}^2}\mathbf{b}_{ks}^H(\sum_{i=1}^K\sum_{j=1}^{S_i}\bar{\beta}_{ij}^2\mathbf{H}_i\mathbf{w}_{ij}\mathbf{w}_{ij}^H\mathbf{H}_i^H + \\ \boldsymbol{\Psi} + \mu_{ks}\mathbf{I}_N)\mathbf{b}_{ks} - \mathbf{b}_{ks}^H\mathbf{H}_k\mathbf{w}_{ks} - \mathbf{w}_{ks}^H\mathbf{H}_k^H\mathbf{b}_{ks} + 1, \ \forall k,s.$$

It implies

$$\mathbf{w}_{ks}^H(\mathbf{H}_k^H\sum_{i=1}^K\sum_{j=1,(i,j)\neq(k,s)}^{S_i}\mathbf{b}_{i,j}\mathbf{b}_{i,j}^H\mathbf{H}_k + \mathbf{R}_{nk})\mathbf{w}_{ks} = \\ \frac{1}{\bar{\beta}_{ks}^2}\mathbf{b}_{ks}^H(\sum_{i=1}^K\sum_{j=1,(i,j)\neq(k,s)}^{S_i}\bar{\beta}_{ij}^2\mathbf{H}_i\mathbf{w}_{ij}\mathbf{w}_{ij}^H\mathbf{H}_i^H + \boldsymbol{\Psi} + \\ \mu_{ks}\mathbf{I}_N)\mathbf{b}_{ks}, \ \forall k,s. \tag{40}$$

Collecting the above expression for all $k$ and $s$ gives

$$(\bar{\mathbf{Y}} + \boldsymbol{\Theta})\bar{\boldsymbol{\beta}}^2 = [a_{11},\cdots,a_{1S_1},\cdots,a_{K1},\cdots,a_{KS_K}]^T = \tilde{\bar{\mathbf{P}}}\mathbf{x} \\ \Rightarrow \bar{\boldsymbol{\beta}}^2 = \boldsymbol{\Theta}^{-1}(\mathbf{I} + \bar{\mathbf{Y}}\boldsymbol{\Theta}^{-1})^{-1}\tilde{\bar{\mathbf{P}}}\mathbf{x} \tag{41}$$

where $\bar{\boldsymbol{\beta}}^2 = [\bar{\beta}_{11}^2,\cdots,\bar{\beta}_{1S_1}^2,\cdots,\bar{\beta}_{K1}^2,\cdots,\bar{\beta}_{KS_K}^2]^T$, $\boldsymbol{\Theta} = \mathrm{diag}(\theta_{11},\cdots,\theta_{1K_1},\cdots,\theta_{K1},\cdots,\theta_{KS_K})$, $a_{ks} = \mathbf{b}_{ks}^H\boldsymbol{\Psi}\mathbf{b}_{ks} + \mu_{ks}\mathbf{b}_{ks}^H\mathbf{b}_{ks}$, $\tilde{\bar{\mathbf{P}}} = [\bar{\bar{\mathbf{P}}},\bar{\mathbf{P}}]$ and $\bar{\mathbf{Y}} = [\bar{\mathbf{y}}_{11},\cdots,\bar{\mathbf{y}}_{1S_1},\cdots,\bar{\mathbf{y}}_{K1}\cdots,\bar{\mathbf{y}}_{KS_K}]^T$ with $\theta_{ks} = \mathbf{w}_{ks}^H\mathbf{R}_{nk}\mathbf{w}_{ks}$, $\bar{\bar{\mathbf{P}}} \in \Re^{S\times N} = |\mathbf{B}^H|^2$, $\bar{\mathbf{P}} = \mathrm{diag}(\bar{p}_{11},\cdots,\bar{p}_{1S_1},\cdots,\bar{p}_{K1},\cdots,\bar{p}_{KS_K})$, $\bar{\mathbf{y}}_{ks} = [-|\mathbf{b}_{ks}^H\mathbf{H}_1\mathbf{w}_{11}|^2,\cdots,\bar{z}_{ks},\cdots,-|\mathbf{b}_{ks}^H\mathbf{H}_K\mathbf{w}_{K1}|,\cdots, -|\mathbf{b}_{ks}^H\mathbf{H}_K\mathbf{w}_{KS_K}|]^T$ and $\bar{z}_{ks} = \mathbf{w}_{ks}^H\mathbf{H}_k^H\sum_{i=1}^K\sum_{j=1,(i,j)\neq(k,s)}^{S_i}\mathbf{b}_{i,j}\mathbf{b}_{i,j}^H\mathbf{H}_k\mathbf{w}_{ks}$. Next we examine two important properties of $(\mathbf{I} + \bar{\mathbf{Y}}\boldsymbol{\Theta}^{-1})^{-1}$. To this end, we examine the following Theorem.

*Theorem 2*: Let $\mathbf{A} \in \Re^{n\times n}$ and $\mathbf{A}_{(i,j),(i\neq j)} \leq 0, 1 \leq i(j) \leq n$. If the diagonal elements of $\mathbf{A}$ are $\mathbf{A}_{(i,i)} = 1 - \sum_{j=1,j\neq i}^n \mathbf{A}_{(j,i)}$, then

$$\text{Property 1}: \mathbf{A}^{-1} \geq 0 \tag{42}$$

$$\text{Property 2}: |||\mathbf{A}^{-1}|||_1 = 1 \tag{43}$$

where $(.) \geq 0$ and $|||.|||_1$ denote matrix non-negativity and one norm, respectively.

*Proof:* See Appendix A. ■

According to the first property of *Theorem 2,* if $\{\theta_{ks} > 0, \forall s\}_{k=1}^K$[5], the inverse of $(\mathbf{I} + \bar{\mathbf{Y}}\boldsymbol{\Theta}^{-1})$ exists and it has nonnegative entries. Consequently, for any positive $\{\psi_n\}_{n=1}^N$ and $\{\mu_{ks}, \forall s\}_{k=1}^K$, $\{\bar{\beta}_{ks}, \forall s\}_{k=1}^K$ of (41) are strictly positive[6]. Now, by selecting $\{\psi_n\}_{n=1}^N$ and $\{\mu_{ks}, \forall s\}_{k=1}^K$ such that (41) is fulfilled, we can transfer the MSE of each symbol from downlink to interference channel ensuring $\{\xi_{ks}^{DL} = \xi_{ks}^{I_1}, \forall s\}_{k=1}^K$, where $\xi_{ks}^{I_1}$ is the MSE of the $k$th user $s$th symbol at step 1 of **Algorithm I**. Here we should also select $\{\psi_n\}_{n=1}^N$ and $\{\mu_{ks}, \forall s\}_{k=1}^K$ such that the power constraint of $\mathcal{P}3$ at step 3 of **Algorithm I** is satisfied. To this end, we examine the steps (2) and (3) of this algorithm.

Like in Section IV, we perform step 2 of **Algorithm 1** by updating $\mathbf{t}_{ks}$ using MMSE receiver as

$$\mathbf{t}_{ks} = (\boldsymbol{\Gamma}_c + \boldsymbol{\Delta}_{ks})^{-1}\mathbf{H}_k\mathbf{v}_{ks}\zeta_{ks} \tag{44}$$

$$= (\sum_{i=1}^K\sum_{j=1}^{S_i}\bar{\beta}_{ij}\mathbf{H}_i\mathbf{w}_{ij}\mathbf{w}_{ij}^H\mathbf{H}_i^H + \boldsymbol{\Psi} + \mu_{ks}\mathbf{I})^{-1}\mathbf{H}_k\mathbf{w}_{ks}\bar{\beta}_{ks}$$

where the second equality is obtained from (39). The above expression shows that by choosing $\{\mu_{ks} > 0, \forall s\}_{k=1}^K$, $\{\psi_n > 0\}_{n=1}^N$, we ensure $(\sum_{i=1}^K\sum_{j=1}^{S_i}\bar{\beta}_{ij}\mathbf{H}_i\mathbf{w}_{ij}\mathbf{w}_{ij}^H\mathbf{H}_i^H + \boldsymbol{\Psi} + \mu_{ks}\mathbf{I})^{-1}$ exists. Next, we transfer the symbol-wise MSE from interference to downlink channel by satisfying the power constraint of $\mathcal{P}3$ (i.e., we perform step 3).

### *B. Symbol-wise MSE transfer (From interference to downlink channel)*

For a given symbol MSE in the interference channel with $\boldsymbol{\zeta} = \mathbf{I}$, we can achieve the same symbol MSE in the downlink channel by using a nonzero scaling factor ($\beta_{ks}$) which satisfies

$$\widetilde{\mathbf{b}}_{ks} = \beta_{ks}\mathbf{t}_{ks}, \quad \widetilde{\mathbf{w}}_{ks} = \mathbf{v}_{ks}/\beta_{ks}. \tag{45}$$

[5] For $\mathcal{P}3$, $\{\mathbf{w}_{ks}^H\mathbf{R}_{nk}\mathbf{w}_{ks} > 0, \forall s\}_{k=1}^K$ is always true.

[6] Note that the application of (43) will be clear in the sequel (see (55)).

Here we use the notations $\widetilde{\mathbf{B}}$ and $\widetilde{\mathbf{W}}$ to differentiate with the precoder and decoder matrices used in Section VI-A. By substituting (45) into $\xi_{ks}^{DL}$ (with $\widetilde{\mathbf{B}}$=$\mathbf{B}$, $\widetilde{\mathbf{W}}$=$\mathbf{W}$), then equating the resulting symbol MSE with that of the interference channel (7) and after some straightforward steps, we get

$$\frac{1}{\beta_{ks}^2}\mathbf{v}_{ks}^H(\mathbf{H}_k^H\sum_{i=1}^{K}\sum_{j=1,(i,j)\neq(k,s)}^{S_i}\beta_{ij}^2\mathbf{t}_{ij}\mathbf{t}_{ij}^H\mathbf{H}_k+\mathbf{R}_{nk})\mathbf{v}_{ks}=$$
$$\mathbf{t}_{ks}^H(\sum_{i=1}^{K}\sum_{j=1,(i,j)\neq(k,s)}^{S_i}\mathbf{H}_i\mathbf{v}_{ij}\mathbf{v}_{ij}^H\mathbf{H}_i^H+\mathbf{\Psi}+\mu_{ks}\mathbf{I})\mathbf{t}_{ks},\forall k,s.$$

By collecting the above equalities for all $k$ and $s$, $\{\beta_{ks},\forall s\}_{k=1}^K$ can be determined by

$$\begin{aligned}(\check{\mathbf{Y}}+\check{\mathbf{\Omega}})\boldsymbol{\beta}^2=&[\mathbf{v}_{11}^H\mathbf{R}_{n1}\mathbf{v}_{11},\cdots,\mathbf{v}_{1S_1}^H\mathbf{R}_{n1}\mathbf{v}_{1S_1},\\&\cdots,\mathbf{v}_{K1}^H\mathbf{R}_{nK}\mathbf{v}_{K1},\cdots,\mathbf{v}_{KS_K}^H\mathbf{R}_{nK}\mathbf{v}_{KS_K}]^T\\=&\mathbf{\Theta}\bar{\boldsymbol{\beta}}^2=\mathbf{\Theta}\mathbf{\Theta}^{-1}(\mathbf{I}+\bar{\mathbf{Y}}\mathbf{\Theta}^{-1})^{-1}\bar{\bar{\mathbf{P}}}\mathbf{x}\\\Rightarrow\boldsymbol{\beta}^2=&(\check{\mathbf{Y}}+\check{\mathbf{\Omega}})^{-1}(\mathbf{I}+\bar{\mathbf{Y}}\mathbf{\Theta}^{-1})^{-1}\bar{\bar{\mathbf{P}}}\mathbf{x}\\=&\check{\mathbf{\Omega}}^{-1}(\mathbf{I}+\check{\mathbf{Y}}\check{\mathbf{\Omega}}^{-1})^{-1}(\mathbf{I}+\bar{\mathbf{Y}}\mathbf{\Theta}^{-1})^{-1}\bar{\bar{\mathbf{P}}}\mathbf{x}\end{aligned}\tag{46}$$

where the third equality follows from (41), $\boldsymbol{\beta}^2=[\beta_{11}^2,\cdots,\beta_{1S_1}^2,\cdots,\beta_{K1}^2,\cdots,\beta_{KS_K}^2]^T$, $\mathbf{\Omega}=\mathrm{diag}(\mathbf{t}_{11}^H\mathbf{\Psi}\mathbf{t}_{11},\cdots,\mathbf{t}_{1S_1}^H\mathbf{\Psi}\mathbf{t}_{1S_1},\cdots,\mathbf{t}_{K1}^H\mathbf{\Psi}\mathbf{t}_{K1},\cdots,\mathbf{t}_{KS_K}^H\mathbf{\Psi}\mathbf{t}_{KS_K})$, $\bar{\mathbf{\Omega}}=\mathrm{diag}(\mu_{11}\mathbf{t}_{11}^H\mathbf{t}_{11},\cdots,\mu_{1S_1}\mathbf{t}_{1S_1}^H\mathbf{t}_{1S_1},\cdots,\mu_{K1}\mathbf{t}_{K1}^H\mathbf{t}_{K1},\cdots,\mu_{KS_K}\mathbf{t}_{KS_K}^H\mathbf{t}_{KS_K})$, $\check{\mathbf{\Omega}}=\mathbf{\Omega}+\bar{\mathbf{\Omega}}$ and $\check{\mathbf{Y}}=[\check{\mathbf{y}}_{11},\cdots,\check{\mathbf{y}}_{1S_1},\cdots,\check{\mathbf{y}}_{K1}\cdots,\check{\mathbf{y}}_{KS_K}]^T$ with $\check{\mathbf{y}}_{ks}=[-|\mathbf{t}_{11}^H\mathbf{H}_1\mathbf{v}_{ks}|^2,\cdots,\check{z}_{ks},\cdots,-|\mathbf{t}_{K1}^H\mathbf{H}_K\mathbf{v}_{ks}|^2,\cdots,-|\mathbf{t}_{KS_K}^H\mathbf{H}_K\mathbf{v}_{ks}|^2]^T$ and $\check{z}_{ks}=\mathbf{t}_{ks}^H\sum_{i=1}^K\sum_{j=1,(i,j)\neq(k,s)}^{S_i}\mathbf{H}_i\mathbf{v}_{ij}\mathbf{v}_{ij}^H\mathbf{H}_i^H\mathbf{t}_{ks}$. By applying *Theorem 2*, it can be shown that $\{\beta_{ks},\forall s\}_{k=1}^K$ are strictly positive for $\{\psi_n>0\}_{n=1}^N$ and $\{\mu_{ks}>0,\forall s\}_{k=1}^K$. The power constraints of the $n$th BS antenna and $k$th user $s$th symbol are given by

$$\tilde{\check{\mathbf{b}}}_n^H\tilde{\check{\mathbf{b}}}_n=\check{\mathbf{t}}_n^H\mathbf{\Upsilon}\check{\mathbf{t}}_n\leq\check{p}_n,\forall n\tag{47}$$
$$\widetilde{\mathbf{b}}_{ks}^H\widetilde{\mathbf{b}}_{ks}=\beta_{ks}^2\mathbf{t}_{ks}^H\mathbf{t}_{ks}\leq\check{p}_{ks},\ \forall k,s\tag{48}$$

where $\mathbf{\Upsilon}=\mathrm{diag}(\beta_{11}^2,\cdots,\beta_{1S_1},\cdots,\beta_{K1}^2,\cdots,\beta_{KS_K})$. Multiplying both sides of (47) by $\psi_n$ and stacking the resulting inequality for all $n$ yields

$$\check{\mathbf{P}}\boldsymbol{\psi}\geq\tilde{\mathbf{\Omega}}\boldsymbol{\beta}^2\tag{49}$$

where $\check{\mathbf{P}}=\mathrm{diag}(\check{p}_1,\cdots,\check{p}_N)$ and $\tilde{\mathbf{\Omega}}=\mathbf{\Psi}|\mathbf{T}|^2$. Like in the above expression, by multiplying both sides of (48) with $\mu_{ks}$ and collecting the resulting inequality for all $k$ and $s$, the power constraints (48) can be expressed as

$$\bar{\check{\mathbf{P}}}\boldsymbol{\mu}\geq\bar{\mathbf{\Omega}}\boldsymbol{\beta}^2\tag{50}$$

where $\bar{\check{\mathbf{P}}}=\mathrm{diag}(\check{p}_{11},\cdots,\check{p}_{1S_1},\cdots,\check{p}_{K1},\cdots,\check{p}_{KS_K})$. By employing $\boldsymbol{\beta}^2$ of (46), (49) and (50) can be combined as

$$\begin{aligned}\mathbf{x}'\geq\bar{\tilde{\mathbf{\Omega}}}\boldsymbol{\beta}^2=&\bar{\tilde{\mathbf{\Omega}}}\check{\mathbf{\Omega}}^{-1}(\mathbf{I}+\check{\mathbf{Y}}\check{\mathbf{\Omega}}^{-1})^{-1}(\mathbf{I}+\bar{\mathbf{Y}}\mathbf{\Theta}^{-1})^{-1}\bar{\bar{\mathbf{P}}}(\bar{\check{\mathbf{P}}})^{-1}\mathbf{x}'\\=&\gimel(\mathbf{x}')\mathbf{x}'\end{aligned}\tag{51}$$

where $\bar{\check{\check{\mathbf{P}}}}=\mathrm{blkdiag}(\check{\mathbf{P}},\bar{\check{\mathbf{P}}})$, $\bar{\tilde{\mathbf{\Omega}}}=[\tilde{\mathbf{\Omega}}^T,\bar{\mathbf{\Omega}}^T]^T$, $\mathbf{x}'=\bar{\check{\check{\mathbf{P}}}}[\boldsymbol{\psi}\ \boldsymbol{\mu}]^T$ and $\gimel(\mathbf{x}')=\bar{\tilde{\mathbf{\Omega}}}\check{\mathbf{\Omega}}^{-1}(\mathbf{I}+\check{\mathbf{Y}}\check{\mathbf{\Omega}}^{-1})^{-1}(\mathbf{I}+\bar{\mathbf{Y}}\mathbf{\Theta}^{-1})^{-1}\bar{\bar{\mathbf{P}}}(\bar{\check{\check{\mathbf{P}}}})^{-1}$.

Next we show that there exists $\mathbf{x}'>0$ such that (51) is satisfied. Towards this end, we consider the following discrete-time switched system [16].

$$\bar{\mathbf{x}}_{n+1}=\mathbf{F}_{\sigma_n}\bar{\mathbf{x}}_n\ \text{for}\ n=0,1,2,\cdots\tag{52}$$

where $\bar{\mathbf{x}}\in\Re^{m\times1}$ is a state, $\mathbf{F}_{\sigma_n}\in\Re^{m\times m}$ is a switching matrix and $\sigma_n\in\{0,1,2,\cdots\}$. According to [16] (Remark 2 of [16]), the above system is marginally stable (convergent) if

$$\max_{\sigma_n}\|\mathbf{F}_{\sigma_n}\|_\star=1\ \text{for}\ n=0,1,\cdots\tag{53}$$

where $\|.\|_\star$ denotes an induced matrix norm.

Let us consider the following iteration

$$\mathbf{x}'_{n+1}=\gimel(\mathbf{x}'_n)\mathbf{x}'_n,\quad\text{for}\ n=0,1,2,\cdots.\tag{54}$$

Now if we assume $\gimel(\mathbf{x}'_n)=\mathbf{F}_{\sigma_n},\forall n$[7], we can interpret (54) as a discrete time switched system. Consequently, the above iteration is guaranteed to converge if $\max_n\|\gimel(\mathbf{x}'_n)\|_\star=1$. It is known that $|||.|||_1$ is an induced matrix norm [17]. For any $\mathbf{x}'$, the matrix one norm of $\gimel(\mathbf{x}')$ is given by

$$\begin{aligned}&|||\gimel(\mathbf{x}')|||_1\\=&|||\bar{\tilde{\mathbf{\Omega}}}\check{\mathbf{\Omega}}^{-1}(\mathbf{I}+\check{\mathbf{Y}}\check{\mathbf{\Omega}}^{-1})^{-1}(\mathbf{I}+\bar{\mathbf{Y}}\mathbf{\Theta}^{-1})^{-1}\bar{\bar{\mathbf{P}}}(\bar{\check{\check{\mathbf{P}}}})^{-1}|||_1\\\leq&|||\check{\check{\mathbf{\Omega}}}|||_1|||(\mathbf{I}+\check{\mathbf{Y}}\check{\mathbf{\Omega}}^{-1})^{-1}|||_1|||(\mathbf{I}+\bar{\mathbf{Y}}\mathbf{\Theta}^{-1})^{-1}|||_1|||\check{\check{\mathbf{P}}}|||_1\\=&|||\check{\check{\mathbf{\Omega}}}|||_1|||\check{\check{\mathbf{P}}}|||_1\leq1\end{aligned}\tag{55}$$

where $\check{\check{\mathbf{\Omega}}}=[\bar{\tilde{\mathbf{\Omega}}}\check{\mathbf{\Omega}}^{-1}\ \mathbf{0}_{(N+S)\times N}]$, $\check{\check{\mathbf{P}}}=[\bar{\bar{\mathbf{P}}}(\bar{\check{\check{\mathbf{P}}}})^{-1};\ \mathbf{0}_{N\times(N+S)}]$, the second inequality is due to the fact that $|||\mathbf{XY}|||_1\leq|||\mathbf{X}|||_1|||\mathbf{Y}|||_1$ [17] (page 290), the third equality is obtained by applying *Theorem 2* and the last inequality employs the following facts. Using the definition (73) (see Appendix A), one can get $|||\check{\check{\mathbf{\Omega}}}|||_1\leq1$ by applying (46) and (51), and $|||\check{\check{\mathbf{P}}}|||_1\leq1$ by applying (13), (41) and (51).

Thus, $\max_n\|\gimel(\mathbf{x}'_n)\|_\star=1$ holds true and (54) is guaranteed to converge. As we can see (54) is derived by using (41) and (46). Thus, the solution of (54) also satisfies (41) and (46). Moreover, for any initial $\mathbf{x}'_0>0$, since $\gimel(\mathbf{x}'_n),\forall n$ is positive, the solution of (54) is strictly positive and $[\boldsymbol{\psi}\ \boldsymbol{\mu}]^T=(\bar{\check{\check{\mathbf{P}}}})^{-1}\mathbf{x}'>0$ which is the desired result.

Once the feasible $\{\mu_{ks},\forall s\}_{k=1}^K$ and $\{\lambda_n\}_{n=1}^N$ are obtained, step 4 of **Algorithm 1** is immediate. As a result, $\mathcal{P}3$ can be solved using **Algorithm I** with an additional power allocation step which will be detailed in Section VIII.

### C. Extension of the current duality for $\mathcal{P}3$ with a total BS power constraint

In this subsection, we show the extension of the current duality for $\mathcal{P}3$ with a total BS power constraint. For this problem, we set $\Delta_{ks}$ of (39) as $\Delta_{ks}=\mathbf{I},\forall k,s$ (i.e., like in Section IV-B). Upon doing so, $\boldsymbol{\beta}^2$ is computed directly from the first equality of (46) (i.e., the bound (55) is not needed). By summing the left and right hand sides of this equality, one can get $\mathrm{tr}\{\widetilde{\mathbf{B}}\widetilde{\mathbf{B}}^H\}=P_{max}$. This shows that for $\mathcal{P}3$ with a total BS power constraint problem, the total BS power at step 3 of **Algorithm I** is satisfied. Thus, one can apply **Algorithm I** (with the additional power allocation step) to solve the latter problem by setting $\Delta_{ks}$ of (39) as $\Delta_{ks} = \mathbf{I}, \forall k, s$.

[7]Since $\gimel(\mathbf{x}'_n)$ is the products of stochastic matrices (see the proof of *Theorem 2*), $\gimel(\mathbf{x}'_n)$ is a bounded matrix for any $\mathbf{x}'>0$. Thus, the assumption $\gimel(\mathbf{x}'_n)=\mathbf{F}_{\sigma_n},\forall n$ holds true.

For other total BS power constrained WMSE-based problems, the current duality based algorithm can be applied like in this subsection. The details are omitted for conciseness. Note that for such problem types, the duality algorithm of the current paper has the same complexity as that of [5].

## VII. User-wise MSE downlink-interference duality

This section establishes user-wise MSE duality between downlink and interference channels. This duality is established to solve the problems of type $\mathcal{P}4$.

### *A. User-wise MSE transfer (From downlink to interference channel)*

To apply this MSE transfer for $\mathcal{P}4$, we set the interference channel precoder, decoder, noise covariance, input covariance and MSE weight matrices as

$$\mathbf{V}_k = \tilde{\beta}_k \mathbf{W}_k,\ \mathbf{T}_k = \mathbf{B}_k/\tilde{\beta}_k, \boldsymbol{\zeta} = \mathbf{I},\ \boldsymbol{\Delta}_{ks} = \boldsymbol{\Psi} + \mu_k \mathbf{I}. \quad (56)$$

Like in Section VI, substituting (56) into (8), equating $\{\xi_k^{DL} = \xi_k^{I}\}_{k=1}^K$ and after some straightforward steps, we get the following system of equations

$$(\tilde{\mathbf{Y}} + \tilde{\boldsymbol{\Theta}})\tilde{\boldsymbol{\beta}}^2 = \tilde{\bar{\mathbf{P}}}\mathbf{x},\ \Rightarrow \tilde{\boldsymbol{\beta}}^2 = \tilde{\boldsymbol{\Theta}}^{-1}(\mathbf{I} + \tilde{\mathbf{Y}}\tilde{\boldsymbol{\Theta}}^{-1})^{-1}\tilde{\bar{\mathbf{P}}}\mathbf{x} \quad (57)$$

where

$$\tilde{\mathbf{Y}}_{(k,l)} = \begin{cases} \sum_{i=1, i\neq k}^K \|\mathbf{W}_k^H \mathbf{H}_k^H \mathbf{B}_i\|_F^2 & \text{for } k = l \\ -\|\mathbf{W}_l^H \mathbf{H}_l^H \mathbf{B}_k\|_F^2, & \text{for } k \neq l \end{cases} \quad (58)$$

$\tilde{\boldsymbol{\Theta}} = \mathrm{diag}(\theta_1, \cdots, \theta_K)$, $\tilde{\boldsymbol{\beta}}^2 = [\tilde{\beta}_1^2, \cdots, \tilde{\beta}_K^2]^T$, $\tilde{\bar{\mathbf{P}}} = [\bar{\bar{\mathbf{P}}}, \tilde{\mathbf{P}}]$ with $\tilde{\theta}_k = \mathrm{tr}\{\mathbf{W}_k^H \mathbf{R}_{nk} \mathbf{W}_k\}$, $\bar{\bar{\mathbf{P}}} \in \Re^{S \times N} = |\mathbf{B}^H|^2$, $\tilde{\mathbf{P}} = \mathrm{diag}(\bar{p}_1, \cdots, \bar{p}_K)$. By applying *Theorem 2*, it can be shown that $\{\tilde{\beta}_k\}_{k=1}^K$ of (57) are strictly positive. Thus, step 1 of **Algorithm I** can be performed using (57). We perform step 2 of **Algorithm I** by updating $\mathbf{t}_{ks}$ using MMSE receiver as

$$\mathbf{t}_{ks} = (\sum_{i=1}^K \tilde{\beta}_i \mathbf{H}_i \mathbf{W}_i \mathbf{W}_i^H \mathbf{H}_i^H + \boldsymbol{\Psi} + \mu_k \mathbf{I})^{-1} \mathbf{H}_k \mathbf{w}_{ks} \tilde{\beta}_k. \quad (59)$$

### *B. User-wise MSE transfer (From interference to downlink channel)*

For a given user MSE in the interference channel with $\boldsymbol{\zeta} = \mathbf{I}$, we can achieve the same MSE in the downlink channel by using nonzero scaling factors ($\{\tilde{\bar{\beta}}_k\}_{k=1}^K$) that satisfy

$$\widetilde{\mathbf{B}}_k = \tilde{\bar{\beta}}_k \mathbf{T}_k, \quad \widetilde{\mathbf{W}}_k = \mathbf{V}_k/\tilde{\bar{\beta}}_k. \quad (60)$$

Here we also use the notations $\widetilde{\mathbf{B}}$ and $\widetilde{\mathbf{W}}$ to differentiate with the precoder and decoder matrices used in Section VII-A. By substituting (60) into $\xi_k^{DL}$ (with $\widetilde{\mathbf{B}}$=$\mathbf{B}$, $\widetilde{\mathbf{W}}$=$\mathbf{W}$) and then equating the resulting user-wise MSE with that of the interference channel (7) and after some steps, $\{\tilde{\bar{\beta}}_k\}_{k=1}^K$ are determined as

$$(\mathbf{I} + \check{\tilde{\mathbf{Y}}}\check{\tilde{\boldsymbol{\Omega}}}^{-1})\check{\tilde{\boldsymbol{\Omega}}}\tilde{\bar{\boldsymbol{\beta}}}^2 = [\mathrm{tr}\{\mathbf{V}_1^H \mathbf{R}_{n1} \mathbf{V}_1\}, \cdots, \mathrm{tr}\{\mathbf{V}_K^H \mathbf{R}_{nK} \mathbf{V}_K\}]^T$$
$$\Rightarrow \tilde{\bar{\boldsymbol{\beta}}}^2 = \check{\tilde{\boldsymbol{\Omega}}}^{-1}(\mathbf{I} + \check{\tilde{\mathbf{Y}}}\check{\tilde{\boldsymbol{\Omega}}}^{-1})^{-1}(\mathbf{I} + \tilde{\mathbf{Y}}\tilde{\boldsymbol{\Theta}}^{-1})^{-1}\tilde{\bar{\mathbf{P}}}\tilde{\mathbf{x}} \quad (61)$$

where the second equality follows from (57), $\tilde{\bar{\boldsymbol{\beta}}}^2 = [\tilde{\bar{\beta}}_1^2, \cdots, \tilde{\bar{\beta}}_K^2]^T$, $\boldsymbol{\Omega}' = \mathrm{diag}(\mathrm{tr}\{\mathbf{T}_1^H \boldsymbol{\Psi} \mathbf{T}_1\}, \cdots, \mathrm{tr}\{\mathbf{T}_K^H \boldsymbol{\Psi} \mathbf{T}_K\})$, $\hat{\boldsymbol{\Omega}} = \mathrm{diag}(\mu_1 \mathrm{tr}\{\mathbf{T}_1^H \mathbf{T}_1\}, \cdots, \mu_K \mathrm{tr}\{\mathbf{T}_K^H \mathbf{T}_K\})$, $\check{\tilde{\boldsymbol{\Omega}}} = \boldsymbol{\Omega}' + \hat{\boldsymbol{\Omega}}$. By applying *Theorem 2*, it can be shown that $\{\tilde{\bar{\beta}}_k\}_{k=1}^K$ are strictly positive for $\{\psi_n > 0\}_{n=1}^N$ and $\{\mu_k > 0\}_{k=1}^K$. Like in Section VI-B, the power constraint of the $n$th BS antenna and $k$th user can thus be expressed as

$$\tilde{\mathbf{x}}' \geq \hat{\bar{\boldsymbol{\Omega}}}\tilde{\bar{\boldsymbol{\beta}}}^2 = \hat{\bar{\boldsymbol{\Omega}}}\check{\tilde{\boldsymbol{\Omega}}}^{-1}(\mathbf{I} + \check{\tilde{\mathbf{Y}}}\check{\tilde{\boldsymbol{\Omega}}}^{-1})^{-1}(\mathbf{I} + \tilde{\mathbf{Y}}\tilde{\boldsymbol{\Theta}}^{-1})^{-1}\tilde{\bar{\mathbf{P}}}\tilde{\mathbf{x}}$$
$$= \tilde{\gimel}(\tilde{\mathbf{x}}')\tilde{\mathbf{x}}' \quad (62)$$

where $\hat{\mathbf{P}} = \mathrm{diag}(\hat{p}_1, \cdots, \hat{p}_K)$, $\hat{\check{\mathbf{P}}} = \mathrm{blkdiag}(\check{\mathbf{P}}, \hat{\check{\mathbf{P}}})$, $\hat{\bar{\boldsymbol{\Omega}}} = [\tilde{\boldsymbol{\Omega}}^T, \hat{\boldsymbol{\Omega}}^T]^T$, $\tilde{\mathbf{x}}' = \check{\tilde{\mathbf{P}}}[\psi\ \ \tilde{\boldsymbol{\mu}}]^T$ and $\tilde{\gimel}(\mathbf{x}') = \hat{\bar{\boldsymbol{\Omega}}}\check{\tilde{\boldsymbol{\Omega}}}^{-1}(\mathbf{I} + \check{\tilde{\mathbf{Y}}}\check{\tilde{\boldsymbol{\Omega}}}^{-1})^{-1}(\mathbf{I} + \tilde{\mathbf{Y}}\tilde{\boldsymbol{\Theta}}^{-1})^{-1}\tilde{\bar{\mathbf{P}}}(\hat{\check{\mathbf{P}}})^{-1}$. Like in Section VI-B, it can be shown that there exists a feasible $\tilde{\mathbf{x}}' > 0$ that satisfy (62) and can be obtained iteratively by

$$\tilde{\mathbf{x}}'_{n+1} = \tilde{\gimel}(\tilde{\mathbf{x}}'_n)\tilde{\mathbf{x}}'_n, \quad \text{for } n = 0, 1, 2, \cdots. \quad (63)$$

By initializing $\mathbf{x}'_0 > 0$, the solution of the above iteration is always positive. Consequently, $\{\lambda_n > 0\}_{n=1}^N$ and $\{\mu_k > 0\}_{k=1}^K$ holds true. Once the feasible $\{\mu_k\}_{k=1}^K$ and $\{\lambda_n\}_{n=1}^N$ are obtained, step 4 of **Algorithm 1** is straightforward. As a result, $\mathcal{P}4$ can be solved using **Algorithm I** with the additional power allocation step of Section VIII.

## VIII. Generalized and improved version of **Algorithm I**

From the discussions of Sections IV - VII, one can understand that each iteration of **Algorithm I** gives a non increasing sequence of symbol (user) WMSE/WSMSE. As can be seen from Section III, the objective of $\mathcal{P}1(\mathcal{P}2)$ is just to minimize the total WSMSE of all symbols (users), whereas the objective of $\mathcal{P}3(\mathcal{P}4)$ is to simultaneously minimize and balance the WMSE of all symbols (users). Thus, **Algorithm I** is appropriate to solve $\mathcal{P}1(\mathcal{P}2)$ of the current paper. For $\mathcal{P}3(\mathcal{P}4)$, although each iteration of **Algorithm I** is able to provide a non increasing sequence of symbol (user) WMSE (i.e., minimizes the maximum WMSE of all symbols (users)), each iteration of this algorithm is not able to guarantee balanced WMSEs of all symbols (users). On the other hand, for an MSE constrained total BS power minimization problem (for example $\mathcal{P}7$ in Section IX), an iterative algorithm that can provide a non increasing sequence of total BS power is required. This shows that **Algorithm I** also can not solve the latter problem. In the following we address the drawbacks of **Algorithm I** just by including a power allocation step into **Algorithm I** as explained below.

In [8], for fixed transmit and receive filters, the power allocation parts of total BS power constrained MSE-based problems have been formulated as GPs by employing the approach and system model of [1] under the assumption that all symbols are strictly active[8]. For this assumption, in [8], we show that the system model of [1] is appropriate to solve any kind of total BS power constrained MSE-based problems using duality approach (alternating optimization). This motivates us to utilize the system model of [1] in the downlink channel only and then include the power allocation step (i.e., GP) into **Algorithm I**. Towards this end, we decompose the precoders and decoders of the downlink channel as

$$\mathbf{B}_k = \mathbf{G}_k \mathbf{P}_k^{1/2}, \quad \mathbf{W}_k = \mathbf{U}_k \boldsymbol{\alpha}_k \mathbf{P}_k^{-1/2}, \ \forall k \tag{64}$$

where $\mathbf{P}_k = \text{diag}(p_{k1}, \cdots, p_{kS_k}) \in \Re^{S_k \times S_k}$, $\mathbf{G}_k = [\mathbf{g}_{k1} \ \cdots \ \mathbf{g}_{kS_k}] \in \mathbb{C}^{N \times S_k}$, $\mathbf{U}_k = [\mathbf{u}_{k1} \ \cdots \ \mathbf{u}_{kS_k}] \in \mathbb{C}^{M_k \times S_k}$ and $\boldsymbol{\alpha}_k = \text{diag}(\alpha_{k1}, \cdots, \alpha_{kS_k}) \in \Re^{S_k \times S_k}$ are the transmit power, unity norm transmit filter, unity norm receive filter and receiver scaling factor matrices of the $k$th user, respectively, i.e., $\{\mathbf{g}_{ks}^H \mathbf{g}_{ks} = \mathbf{u}_{ks}^H \mathbf{u}_{ks} = 1, \forall s\}_{k=1}^K$.

By employing (64) and stacking $\xi = [\xi_{1,1}^{DL}, \cdots, \xi_{K,S_K}^{DL}]^T = [\xi_1^{DL}, \cdots, \xi_S^{DL}]^T = [\{\xi_l^{DL}\}_{l=1}^S]^T$, the $l$th downlink symbol MSE can be expressed as (see [1] and [10] for more details about (64) and the above descriptions)

$$\xi_l^{DL} = p_l^{-1}[(\mathbf{D} + \boldsymbol{\alpha}^2 \boldsymbol{\Phi}^T)\mathbf{p}]_l + p_l^{-1} \alpha_l^2 \mathbf{u}_l^H \mathbf{R}_n \mathbf{u}_l \tag{65}$$

where

$$\boldsymbol{\Phi}_{(l,j)} = \begin{cases} |\mathbf{g}_l^H \mathbf{H} \mathbf{u}_j|^2, & \text{for } l \neq j \\ 0, & \text{for } l = j \end{cases} \tag{66}$$

$$\mathbf{D}_{(l,l)} = \alpha_l^2 |\mathbf{g}_l^H \mathbf{H} \mathbf{u}_l|^2 - 2\alpha_l \Re(\mathbf{u}_l^H \mathbf{H}^H \mathbf{g}_l) + 1, \tag{67}$$

$1 \leq l(j) \leq S$, $\mathbf{P} = \text{blkdiag}(\mathbf{P}_1, \cdots, \mathbf{P}_K) = \text{diag}(p_1, \cdots, p_S)$, $\mathbf{p} = [p_1, \cdots, p_S]^T$, $\mathbf{G} = [\mathbf{G}_1, \cdots, \mathbf{G}_K] = [\mathbf{g}_1, \cdots, \mathbf{g}_S]$, $\mathbf{U} = \text{blkdiag}(\mathbf{U}_1, \cdots, \mathbf{U}_K) = [\mathbf{u}_1, \cdots, \mathbf{u}_S]$ and $\boldsymbol{\alpha} = \text{blkdiag}(\boldsymbol{\alpha}_1, \cdots, \boldsymbol{\alpha}_K) = \text{diag}(\alpha_1, \cdots, \alpha_S)$ with $\|\mathbf{g}_l\|_2 = \|\mathbf{u}_l\|_2 = 1$. Using (65), for fixed $\mathbf{G}, \mathbf{U}$ and $\boldsymbol{\alpha}$, the power allocation part of $\mathcal{P}1$ can be formulated as

$$\min_{\{p_l\}_{l=1}^S} \sum_{l=1}^S \eta_l \xi_l^{DL}, \ \text{s.t} \ \boldsymbol{\varsigma}_n^T \mathbf{p} \leq \breve{p}_n, \ p_l \leq \breve{p}_l \ \forall n, l \tag{68}$$

where $\boldsymbol{\varsigma}_n^T \in \Re^{1 \times S} = \{|[\mathbf{G}_{(n,i)}|^2\}_{i=1}^S$, $[\eta_1, \cdots, \eta_S]^T = [\eta_{11}, \cdots, \eta_{KS_K}]^T$ and $[\breve{p}_l, \cdots, \breve{p}_S]^T = [\breve{p}_{11}, \cdots, \breve{p}_{KS_K}]^T$. As $\xi_l^{DL}$ is a posynomial (where $\{p_l\}_{l=1}^S$ are the variables), (68) is a GP for which global optimality is guaranteed. Thus, it can be efficiently solved using interior point methods with a worst-case polynomial-time complexity [18].

For fixed $\mathbf{G}, \mathbf{U}$ and $\boldsymbol{\alpha}$, the power allocation parts of $\mathcal{P}2 - \mathcal{P}4$ can be formulated as GPs like in $\mathcal{P}1$. Our duality based algorithm for each of these problems including the power allocation step is summarized in **Algorithm II**.

**Algorithm II**

Initialization: Like in **Algorithm I**.

**Repeat**: **Interference channel**

1) For $\mathcal{P}1$ and $\mathcal{P}2$, set $\mathbf{V} = \mathbf{W}, \mathbf{T} = \mathbf{B}$ (i.e., $\bar{\beta} = \tilde{\beta} = 1$), then compute $\{\psi_n, \mu_{ks}, \forall k, s, n\}$ and $\{\psi_n, \mu_k, \forall k, n\}$ using (25) and (38), respectively. For $\mathcal{P}3$ and $\mathcal{P}4$, first compute $\{\psi_n, \mu_{ks}, \forall k, s, n\}$ and $\{\psi_n, \mu_k, \forall k, n\}$ using (54) and (63), respectively, then transfer each symbol and user MSE from downlink to interference channels by (39) and (56), respectively.
2) Update the MMSE receivers of the interference channel for $\mathcal{P}1$, $\mathcal{P}2$, $\mathcal{P}3$ and $\mathcal{P}4$ using (19), (32), (44) and (59), respectively.

**Downlink channel**

3) Transfer the MSE (weighted sum, user or symbol MSE) from interference to downlink channel using (20), (33), (45) and (60) for $\mathcal{P}1$, $\mathcal{P}2$, $\mathcal{P}3$ and $\mathcal{P}4$, respectively.
4) For each of the problems $\mathcal{P}1 - \mathcal{P}4$, decompose the precoder and decoder matrices of each user as in (64). Then, formulate and solve the GP power allocation part. For example, the power allocation part of $\mathcal{P}1$ can be expressed in GP form as (68).
5) For each of the problems $\mathcal{P}1 - \mathcal{P}4$, by keeping $\{\mathbf{P}_k\}_{k=1}^K$ constant, update the receive filters $\{\mathbf{U}_k\}_{k=1}^K$ and scaling factors $\{\boldsymbol{\alpha}_k\}_{k=1}^K$ by applying downlink MMSE receiver approach i.e., $\{\mathbf{U}_k \boldsymbol{\alpha}_k = (\mathbf{H}_k^H \mathbf{G} \mathbf{P} \mathbf{G}^H \mathbf{H}_k + \mathbf{R}_{nk})^{-1} \mathbf{H}_k^H \mathbf{G}_k \mathbf{P}_k\}_{k=1}^K$. Note that in these expressions, $\{\boldsymbol{\alpha}_k\}_{k=1}^K$ are chosen such that each column of $\{\mathbf{U}_k\}_{k=1}^K$ has unity norm. Then, compute $\{\mathbf{B}_k, \mathbf{W}_k\}_{k=1}^K$ by (64).

**Until** convergence.

**Convergence**: It can be shown that at each iteration of this algorithm, the objective function of each of the problems $\mathcal{P}1$ - $\mathcal{P}4$ is non-increasing [4], [7], [19]. Thus, the above iterative algorithm is convergent. However, since $\mathcal{P}1$ - $\mathcal{P}4$ are non-convex, this iterative algorithm is not guaranteed to converge to the global optimum.

In this algorithm, we stop iteration (i.e., our convergence condition) when the difference between the objective functions in two consecutive iterations is smaller than some small value $\tilde{\epsilon}$ (we use $\tilde{\epsilon} = 10^{-6}$ for the simulation).

**Computational complexity**: As can be seen from this algorithm, when we increase the number of users and/or (BS and/or MS antennas), the number and size of optimization variables increase. Because of this, the computational complexity of **Algorithm II** increases as $K$ and/or $N$ and/or $M$ increases. However, studying the complexity of this algorithm as a function of $K, N$ and $M$ needs effort and time. And such a task is beyond the scope of this work and is an open research topic.

The power allocation step of **Algorithm II** has thus the following benefits: (1) For BS power constrained WSMSE minimization problems, this step improves the convergence speed of **Algorithm II** compared to that of **Algorithm I**[9] (for example in $\mathcal{P}1 - P2$). The degree of improvement depends on different parameters (for example $\mathbf{H}_k$, $\boldsymbol{\Delta}_{ks}, \forall k, s$ etc). Thus, the theoretical comparison of these two algorithms in terms of convergence speed requires time and effort. And this task is beyond the scope of this work and it is an open research topic.

[8]Note that this assumption is not always true for all MSE-based problems. However, as mentioned in [1], in practice replacing zero powers by a small value will not affect the overall optimization. Due to this reason, we replace zero powers by $10^{-6}$ in the simulation section.

[9]This is at the expense of additional computation. However, as mentioned in [1] (see Appendix A of [1]), a small desktop computer can solve a GP of 100 variables and 10000 constraints by standard interior point method under a minute. Thus, we believe that the complexity of **Algorithm I** and **Algorithm II** are almost the same.

(2) For symbol-wise (user-wise) WMSE balancing problems, this step helps to balance the WMSE of all symbols (users) (for example in $\mathcal{P}3 - P4$). (3) For MSE constrained total BS power minimization problems, this step ensures a non increasing total BS power at each iteration of **Algorithm II**.

## IX. Application of the proposed duality based algorithm for other problems

### A. MSE based problem with entry-wise power constraint

The symbol-wise WSMSE minimization constrained with entry wise power i.e,. $b_{ksn}^H b_{ksn} \leq \bar{\bar{p}}_{ksn}, \forall k, s, n$ problem is formulated as

$$\mathcal{P}5: \min_{\{\mathbf{B}_k, \mathbf{W}_k\}_{k=1}^K} \sum_{k=1}^{K} \sum_{s=1}^{S_k} \eta_{ks} \xi_{ks}^{DL}, \text{ s.t } b_{ksn}^H b_{ksn} \leq \bar{\bar{p}}_{ksn}. \quad (69)$$

It can be shown that this problem can be solved by **Algorithm II** with $\{\mathbf{\Delta}_{ks} = \text{diag}(\delta_{ks1}, \cdots, \delta_{ksN}), \forall s\}_{k=1}^K$.

### B. Weighted sum rate optimization constrained with per antenna and symbol power problem

By employing the approach of [11] (see (16) of [11]), one can equivalently express the weighted sum rate maximization constrained with per antenna and symbol power problem as

$$\mathcal{P}6: \quad (70)$$

$$\min_{\{\bar{\tau}_{ks}, \bar{\nu}_{ks}, \mathbf{b}_{ks}, \mathbf{w}_{ks}, \forall s\}_{k=1}^K} \sum_{k=1}^{K} \sum_{s=1}^{S_k} \bar{\theta}_{ks} \frac{1}{\bar{\tau}_{ks}} \bar{\nu}_{ks}^{\bar{\gamma}_{ks}} + \sum_{k=1}^{K} \sum_{s=1}^{S_k} \bar{\eta}_{ks} \xi_{ks}^{DL},$$

$$\text{s.t } [\mathbf{B}\mathbf{B}^H]_{(n,n)} \leq \breve{p}_n, \mathbf{b}_{ks}^H \mathbf{b}_{ks} \leq \breve{p}_{ks}, \prod_{k=1}^{K} \prod_{s=1}^{S_k} \bar{\nu}_{ks} = 1, \bar{\tau}_{ks} > 0$$

where $\{0 < \bar{\omega}_{ks} < 1, \forall s\}_{k=1}^K$ are the rate weighting factors for all symbols, $\bar{\eta}_{ks} = \bar{\tau}_{ks}^{\bar{\mu}_{ks}}$, $\bar{\gamma}_{ks} = \frac{1}{1-\bar{\omega}_{ks}}$ $\bar{\mu}_{ks} = \frac{1}{\bar{\omega}_{ks}} - 1$ and $\bar{\theta}_{ks} = \bar{\omega}_{ks} \bar{\mu}_{ks}^{(1-\bar{\omega}_{ks})}$. For fixed $\{\bar{\tau}_{ks}, \bar{\nu}_{ks}, \forall s\}_{k=1}^K$, the above optimization problem has the same mathematical structure as that of $\mathcal{P}1$. Thus, by keeping $\{\bar{\tau}_{ks}, \bar{\nu}_{ks}, \forall s\}_{k=1}^K$ constant, $\{\mathbf{b}_{ks}, \mathbf{w}_{ks}, \forall s\}_{k=1}^K$ can be optimized by applying the MSE duality discussed in Section IV. Moreover, $\{\bar{\tau}_{ks}, \bar{\nu}_{ks}, \forall s\}_{k=1}^K$ and the power allocation part of the above problem can be optimized by a GP method like in (25) of [20]. Consequently, we can apply **Algorithm II** to solve (70). The detailed explanations are omitted for conciseness. The following problems can also be solved by simple modification of **Algorithm II**

$$\begin{aligned}
\mathcal{P}7: & \min_{\{\mathbf{B}_k, \mathbf{W}_k\}_{k=1}^K} \sum_{k=1}^{K} \text{tr}\{\mathbf{B}_k \mathbf{B}_k^H\}, \\
& \text{s.t } [\mathbf{B}\mathbf{B}^H]_{(n,n)} \leq \breve{p}_n, \\
& \text{tr}\{\mathbf{B}_k^H \mathbf{B}_k\} \leq \hat{p}_k \\
& SINR_{ks} \geq \varrho_{ks}, \qquad \forall n, k, s \\
\equiv: & \min_{\{\mathbf{B}_k, \mathbf{W}_k\}_{k=1}^K} \sum_{k=1}^{K} \text{tr}\{\mathbf{B}_k \mathbf{B}_k^H\}, \\
& \text{s.t } [\mathbf{B}\mathbf{B}^H]_{(n,n)} \leq \breve{p}_n, \\
& \text{tr}\{\mathbf{B}_k^H \mathbf{B}_k\} \leq \hat{p}_k \\
& \xi_{ks}^{DL} \leq (1 + \varrho_{ks})^{-1}, \qquad \forall n, k, s
\end{aligned}$$

$$\begin{aligned}
\mathcal{P}8: & \max_{\{\mathbf{B}_k, \mathbf{W}_k\}_{k=1}^K} \min R_{ks} \\
& \text{s.t } [\mathbf{B}\mathbf{B}^H]_{(n,n)} \leq \breve{p}_n, \ \mathbf{b}_{ks}^H \mathbf{b}_{ks} \leq \breve{p}_{ks}, \ \forall n, k, s \\
\equiv: & \min_{\{\mathbf{B}_k, \mathbf{W}_k\}_{k=1}^K} \max \xi_{ks}^{DL} \\
& \text{s.t } [\mathbf{B}\mathbf{B}^H]_{(n,n)} \leq \breve{p}_n, \ \mathbf{b}_{ks}^H \mathbf{b}_{ks} \leq \breve{p}_{ks}, \ \forall n, k, s
\end{aligned}$$

where $SINR_{ks}(R_{ks})$ is the SINR (rate) of the $k$th user $s$th symbol, and we use the fact that $R_{ks} = \log(1 + SINR_{ks})$ and $\xi_{ks}^{DL} = (1 + SINR_{ks})^{-1}$ [1]. It is clearly seen that the application of **Algorithm II** is not limited to the problems of this paper.

Note that under imperfect channel state information (CSI) condition, the stochastic robust design versions of $\mathcal{P}1$ - $\mathcal{P}5$ can be solved like in [10]. However, to the best of our knowledge, the relationship between rate (SINR) and MSE is not known when the CSI is imperfect [7]. Hence, solving the rate (SINR)-based robust design problems (for example, robust versions of $\mathcal{P}6$ - $\mathcal{P}8$) by our duality approach is an open problem.

## X. Simulation Results

In this section, we present simulation results for $\mathcal{P}1 - \mathcal{P}4$. All of our simulation results are averaged over $100$ randomly chosen channel realizations. We set $K = 2$, $N = 4$ and $\{M_k = S_k = 2, \eta_{ks} = \rho_{ks} = \tilde{\eta}_k = \tilde{\rho}_k = 1, \forall s\}_{k=1}^K$. It is assumed that $\mathbf{R}_{n1} = \sigma_1^2 \mathbf{I}_{M_1}$, $\mathbf{R}_{n2} = \sigma_2^2 \mathbf{I}_{M_2}$ and $\sigma_2^2 = 2\sigma_1^2$. The maximum power of each BS antenna is set to $\{\breve{p}_n = 2.5mW\}_{n=1}^N$. And the maximum power allocated to each symbol and user are set to $\{\breve{p}_{ks} = 2.5mW, \forall s\}_{k=1}^K$ and $\{\hat{p}_k = 5mW\}_{k=1}^K$, respectively. For better exposition, we define the Signal-to-noise ratio (SNR) as $P_{\max}/K\sigma_{av}^2$ and it is controlled by varying $\sigma_{av}^2$, where $P_{\max} = 10$mW is the total maximum BS power and $\sigma_{av}^2 = (\sigma_1^2 + \sigma_2^2)/2$. We also compare **Algorithm II** and the algorithm in [2][10].

Note that the algorithm in [2] is designed for coordinated BS systems scenario. And, the iterative algorithm of [2] is based on the per BS power constraint. However, according to [2] and [21], $B$ coordinated BS systems each with $Z$ antennas can be treated as one multiuser MIMO system with $BZ$ antennas. Thus, when $Z = 1$, the considered problem has exactly the same structure as that of [2]. To the best of our knowledge, there is no other general linear algorithm that can solve the problem types $\mathcal{P}1 - \mathcal{P}4$. On the other hand, in all problems, since there are more than one power constraints (i.e., per antenna and symbol (user) powers), all power constraints may not be active at the optimal solution. Due to these reasons, we compare **Algorithm II** and the algorithm in [2] both in terms of the achieved MSE (i.e., minimized MSE) and total utilized BS power at the achieved MSE.

### A. Simulation results for problems $\mathcal{P}1 - \mathcal{P}2$

In this subsection, we compare the performance of our proposed algorithm with that of [2]. As can be seen from Fig.

[10] As will be clear in the sequel, the proposed algorithm and the algorithm of [2] may not utilize the entire 10mW for all noise levels. Thus, the aforementioned SNR can be considered as the maximum SNR of the transmitted signal.

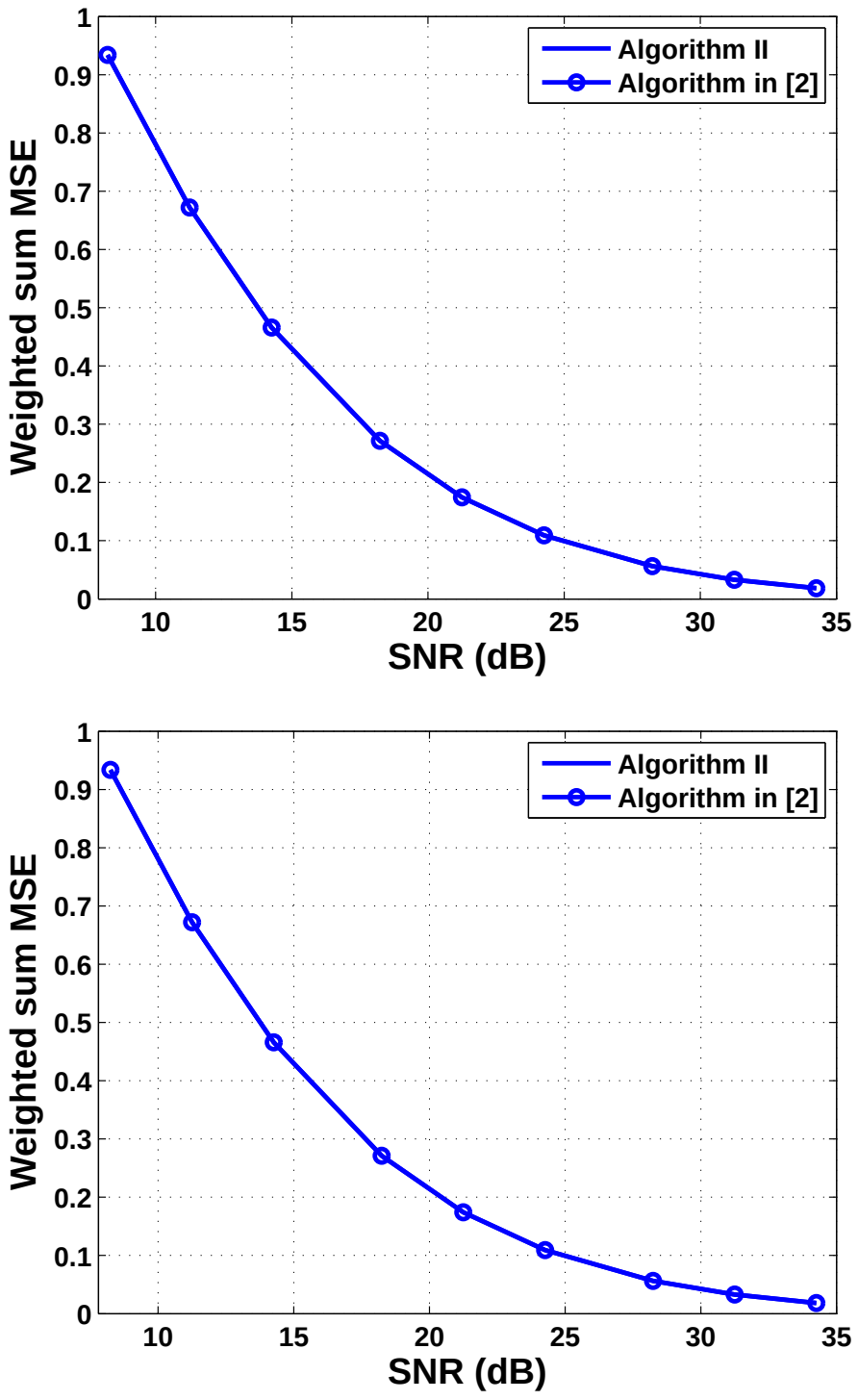

Fig. 2. Comparison of the proposed algorithm (**Algorithm II**) and the algorithm of [2] in terms of WSMSE for: [top] $\mathcal{P}1$, [bottom] $\mathcal{P}2$.

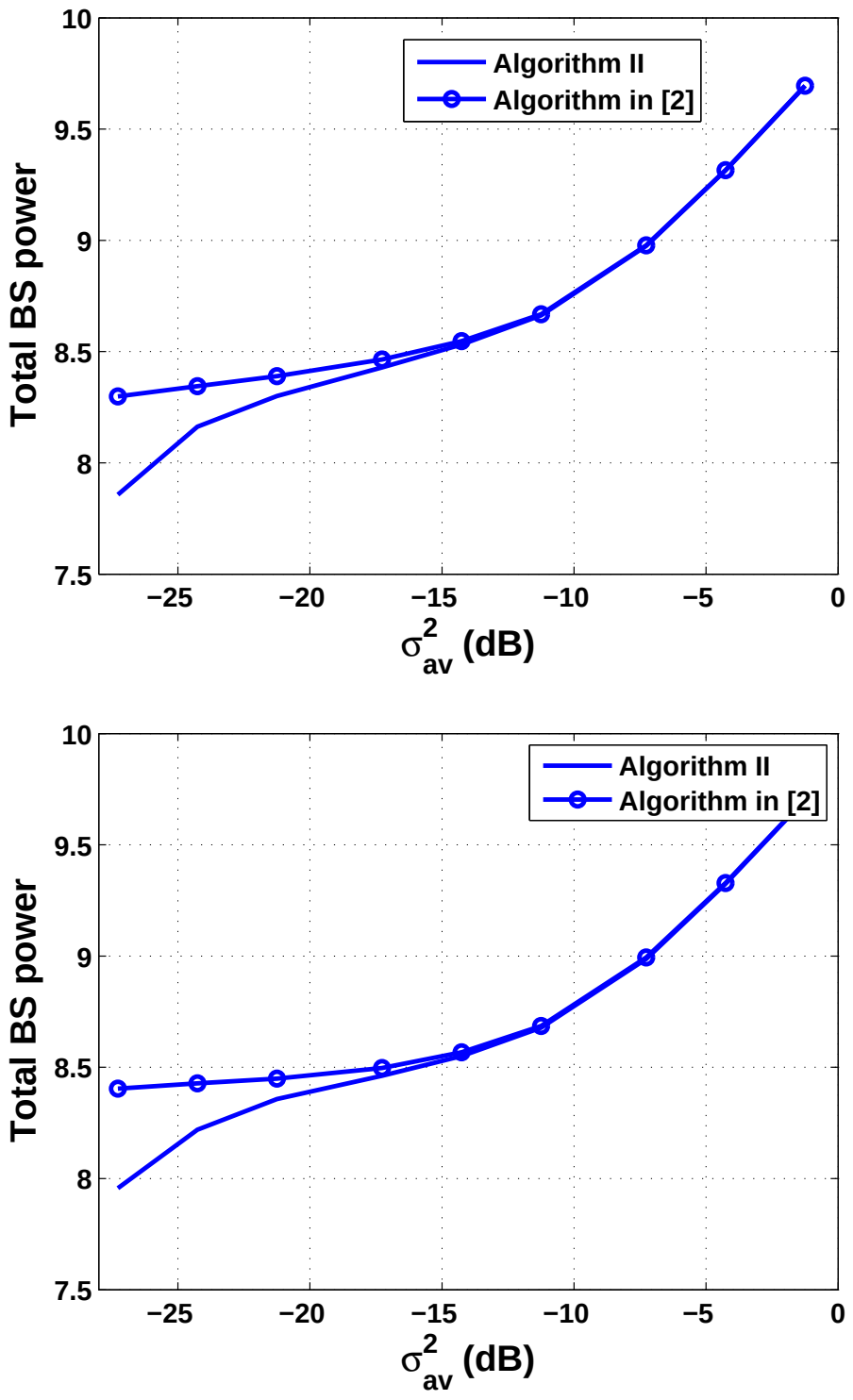

Fig. 3. Comparison of the proposed algorithm (**Algorithm II**) and the algorithm of [2] in terms of total BS power for: [top] $\mathcal{P}1$, [bottom] $\mathcal{P}2$. For this figure we compute $\sigma_{av}^2(dB)$ as $\sigma_{av}^2$(dB)=$10\log\frac{\sigma_{av}^2}{1mW}$.

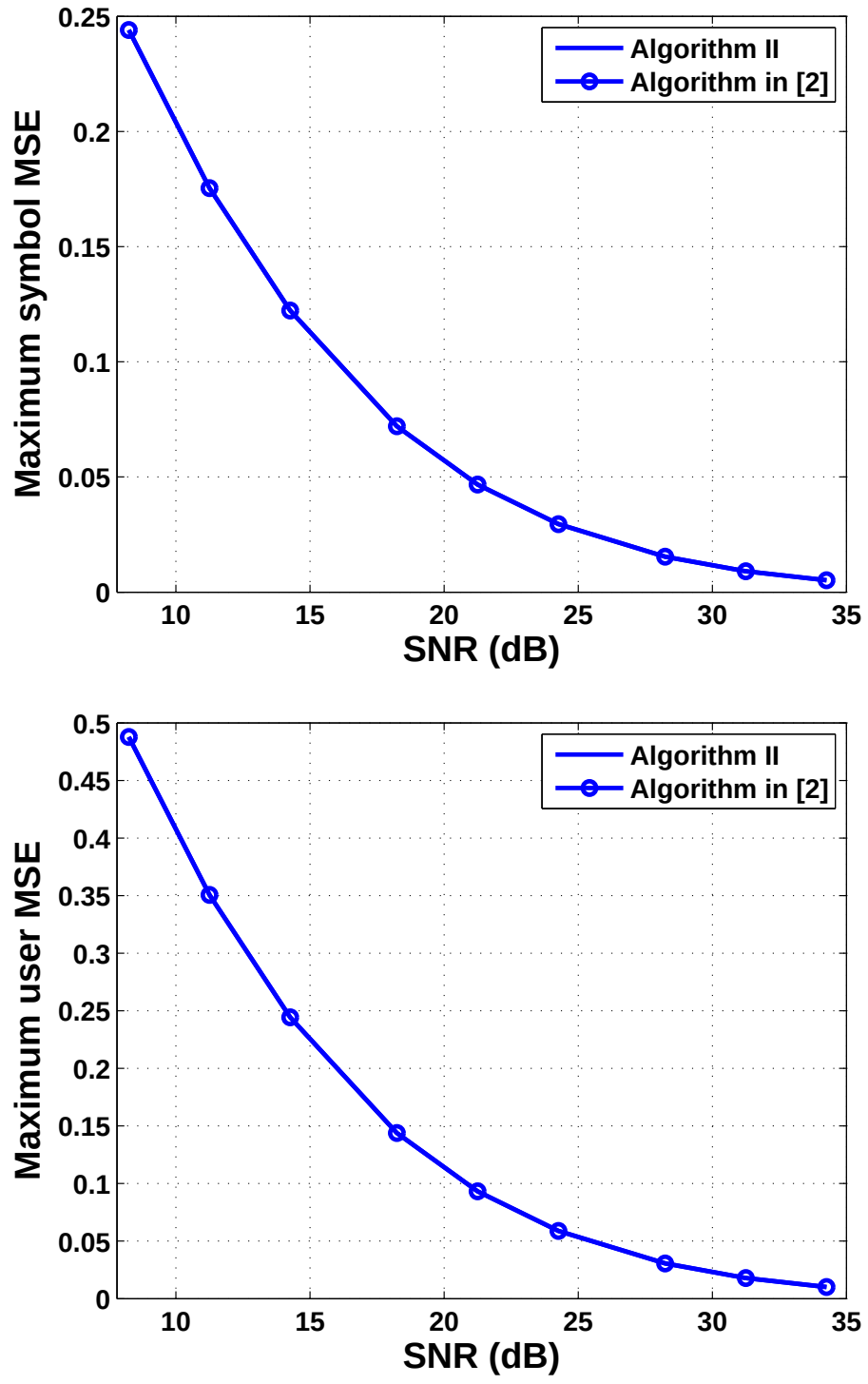

Fig. 4. Comparison of the proposed algorithm (**Algorithm II**) and that of in [2] in terms of maximum achieved MSE for: [top] $\mathcal{P}3$, [bottom] $\mathcal{P}4$.

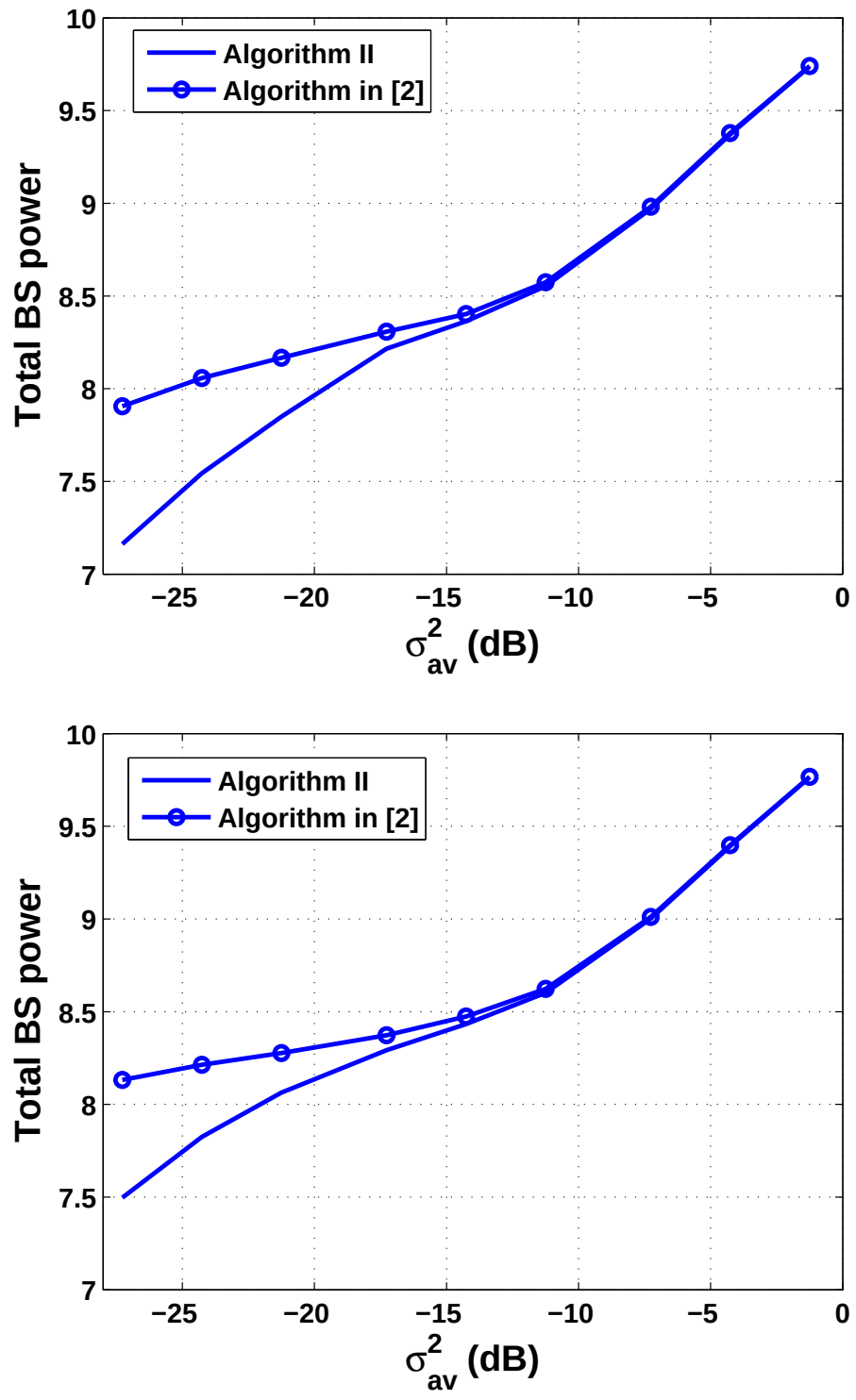

Fig. 5. Comparison of the proposed algorithm (**Algorithm II**) and the algorithm of [2] in terms of total BS power for: [top] $\mathcal{P}3$, [bottom] $\mathcal{P}4$. In this figure we compute $\sigma_{av}^2(dB)$ as $\sigma_{av}^2$(dB)=$10\log\frac{\sigma_{av}^2}{1mW}$.

2, the proposed algorithm and the algorithm in [2] achieve the same symbol-wise and user-wise WSMSEs. Next, we plot the total utilized powers at the BS to achieve these WSMSEs which is shown in Fig. 3. From these two figures, one can see that to achieve the same WSMSE, the proposed duality based iterative algorithm requires less total BS power than that of [2]. This scenario fits to that of [10] and [19] where the sum MSE minimization constrained with a per BS antenna power problem has been examined by duality approach.

### B. Simulation results for problems $\mathcal{P}3 - \mathcal{P}4$

Like in the above subsection, here we compare the performance of our proposed algorithm with that of [2]. For these problems, we also observe from Fig. 4 that the proposed algorithm and the algorithm in [2] achieve the same maximum symbol and user MSEs. And from Fig. 5 the proposed duality based algorithms utilize less total BS power compared to that of [2]. For all of our problems, we observe that to achieve the same MSE, the proposed duality based iterative algorithm utilizes less total BS power compared to the algorithm of [2]. This scenario has also been observed for other MSE and rate based problems in [10], [19], [20].

Note that the problems of [2] are examined directly in the downlink channel. According to [5], [13], in general, duality approach of solving downlink transceiver design problems has easier to handle mathematical structure (lower complexity) compared to that of the direct approach in [2]. For this reason, we believe that the computational complexity of the proposed duality algorithm is not higher than that of [2].

### C. Convergence speed of **Algorithm II**

As can be seen from Section VIII, the overall computational complexity of **Algorithm II** depends on the number of iterations to achieve convergence. In general, the number of iterations to achieve convergence may not be the same for all problems. On the other hand, for each problem, getting the exact number of iterations to achieve convergence analytically is very difficult. Due to these reasons, we provide numerical simulations to demonstrate the convergence speed of **Algorithm II** for $\mathcal{P}1$. As can be seen from Fig. 6, the proposed algorithm converges within few iterations in low, medium and high SNR regions.

## XI. Conclusions

In this paper, we examine different transceiver design problems for multiuser MIMO systems under generalized linear power constraints. The problems are solved for the practically relevant scenario where the noise vector of each MS is a ZMCSCG random variable with arbitrary covariance matrix. For all of our problems, we propose new downlink-interference duality based iterative solutions. The current duality are established by formulating the noise covariance matrices of the dual interference channels as fixed point functions and marginally stable (convergent) discrete-time-switched systems. We show that the proposed duality based iterative algorithms can be extended straightforwardly to solve several practically relevant linear transceiver design problems. We also show that our new MSE downlink-interference duality unify all existing MSE duality. Our simulation results demonstrate that the proposed duality based algorithms utilize less total BS power than that of existing algorithms.

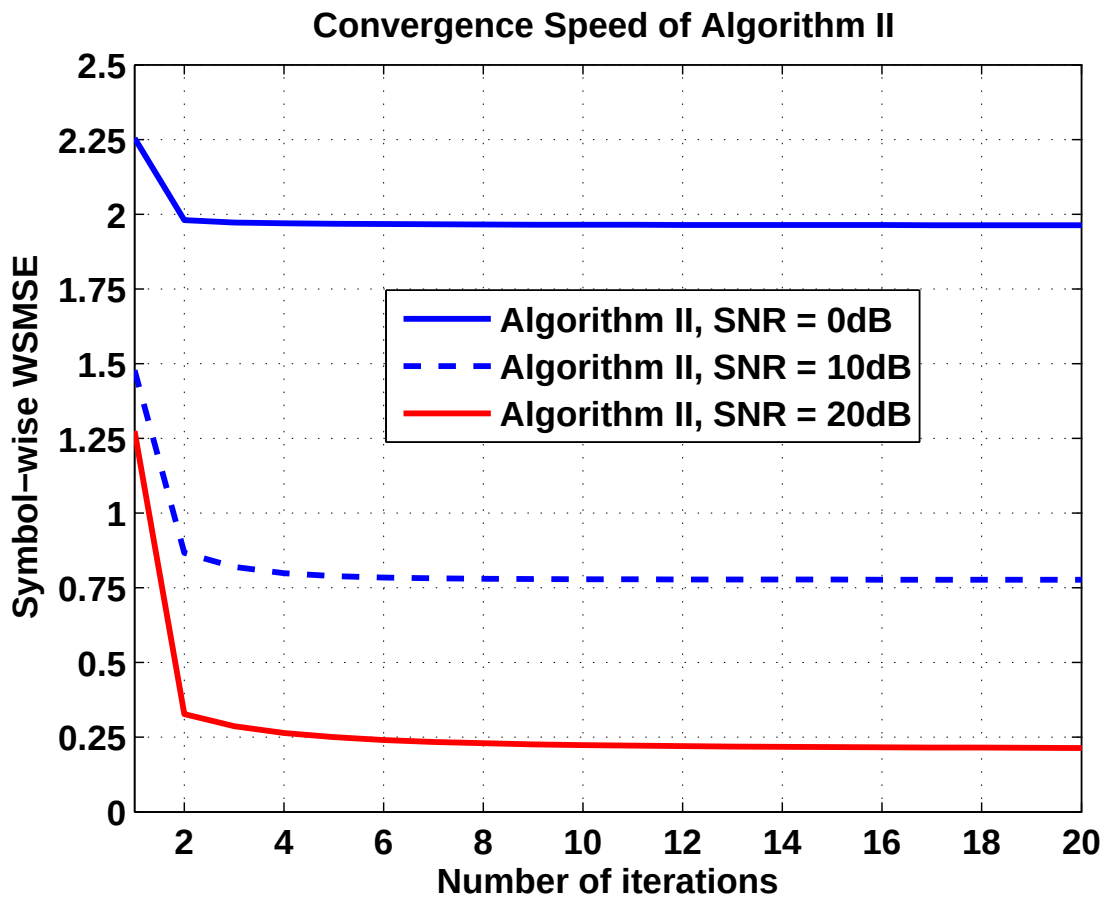

Fig. 6. Convergence speed of **Algorithm II** for $\mathcal{P}1$.

## Appendix A
## Proof of Theorem 2

*Proof:* Define $\tilde{\mathbf{D}} \triangleq \text{diag}(\mathbf{A}_{1,1}, \cdots, \mathbf{A}_{n,n})$ and $\tilde{\mathbf{A}} \triangleq \tilde{\mathbf{D}} - \mathbf{A}$. It follows

$$\mathbf{A} = \tilde{\mathbf{D}} - \tilde{\mathbf{A}} \Rightarrow \mathbf{A}^{-1} = \tilde{\mathbf{D}}^{-1}(\mathbf{I} - \bar{\mathbf{A}})^{-1}$$

where $\bar{\mathbf{A}} = \tilde{\mathbf{A}}\tilde{\mathbf{D}}^{-1}$. Since $(\mathbf{I} - \bar{\mathbf{A}})$ is strictly diagonally dominant matrix, $(\mathbf{I} - \bar{\mathbf{A}})^{-1}$ exists [17] (page 349). Furthermore, if $\rho(\bar{\mathbf{A}}) < 1$, $(\mathbf{I} - \bar{\mathbf{A}})^{-1}$ can be expressed as

$$(\mathbf{I} - \bar{\mathbf{A}})^{-1} = \sum_{k=0}^{\infty} \bar{\mathbf{A}}^k. \tag{71}$$

It follows

$$\bar{\mathbf{A}}^{-1} = \tilde{\mathbf{D}}^{-1}(\mathbf{I} - \bar{\mathbf{A}})^{-1} = \tilde{\mathbf{D}}^{-1}\sum_{k=0}^{\infty} \bar{\mathbf{A}}^k \geq 0 \tag{72}$$

From this equation we can see that if $\rho(\bar{\mathbf{A}}) < 1$, the nonnegativity of $\bar{\mathbf{A}}^{-1}$ can be ensured. Next we show that $\rho(\bar{\mathbf{A}})$ is indeed less than 1. For any $n \times n$ matrix $\mathbf{X}$, we have [17] (pages 294 and 297 of [17])

$$\rho(\mathbf{X}) \leq |||\mathbf{X}|||, \quad |||\mathbf{X}||_1 \triangleq \max_{1\leq j\leq n} \sum_{i=1}^{n} |x_{ij}| \tag{73}$$

where $|||.|||$ is any matrix norm and $|||.|||_1$ is a matrix one norm. By using (73), we get the following bound [17]

$$\rho(\bar{\mathbf{A}}) \leq |||\bar{\mathbf{A}}|||_1 < 1. \tag{74}$$

Since $\mathbf{A}^{-1}$ has nonnegative elements, $\mathbf{A}$ is also an M-matrix [22]. By defining $\mathbf{S} \triangleq \mathbf{A}^{-1}$ and $\mathbf{e} \triangleq \mathbf{1}^{n\times 1}$, we get

$$\begin{aligned}\mathbf{e}^T\mathbf{A} = \mathbf{e}^T \ \Rightarrow \mathbf{e}^T &= \mathbf{e}^T\mathbf{S} = [\sum_{j=1}^{n}\mathbf{S}_{j,1}, \cdots, \sum_{j=1}^{n}\mathbf{S}_{j,n}] \\ \Rightarrow |||\mathbf{S}|||_1 &= 1\end{aligned} \tag{75}$$

where the third equality follows from the fact that $\mathbf{S}$ is a nonnegative matrix. ∎